\documentclass[a4paper,12pt,twoside]{article}
\usepackage{amsmath,amssymb,ifthen}
\usepackage[amsmath,thmmarks]{ntheorem}
\usepackage{paper}
\usepackage[all]{xy}
\usepackage{mathrsfs}

\newcommand{\ad}{\mathrm{ad}}
\DeclareMathOperator{\Spf}{Spf}
\DeclareMathOperator{\Spa}{Spa}
\DeclareMathOperator{\Fix}{Fix}
\DeclareMathOperator{\fix}{fix}

\DeclareMathOperator{\Gys}{Gys}
\DeclareMathOperator{\Adj}{Adj}
\DeclareMathOperator{\loc}{loc}
\DeclareMathOperator{\Nrd}{Nrd}
\DeclareMathOperator{\Sh}{Sh}
\newcommand{\Lotimes}{\overset{\mathbb{L}}{\otimes}}
\newcommand{\Lboxtimes}{\overset{\mathbb{L}}{\boxtimes}}
\newcommand{\M}{\breve{\mathscr{M}}}

\SelectTips{cm}{11}
\title{Lefschetz trace formula for open adic spaces}
\author{Yoichi Mieda}
\begin{document}
\maketitle

\begin{firstfootnote}
 The Hakubi Center for Advanced Research / Department of Mathematics, Kyoto University, Kyoto, 606--8502, Japan

 E-mail address: \texttt{mieda@math.kyoto-u.ac.jp}

\end{firstfootnote}

\begin{abstract}
 In this article, we discuss the Lefschetz trace formula for an adic space
 which is separated smooth of finite type but not necessarily proper over an algebraically closed
 non-archimedean field. 
 Under the condition that there is no set-theoretical fixed point on the boundary,
 we obtain a fixed point formula.
 As an application, we can establish a trace formula for some formal schemes,
 which is applicable to some Rapoport-Zink towers.
 A partial generalization of Fujiwara's trace formula for contracting morphisms is also given.
\end{abstract}

\section{Introduction}
In this paper, we consider the Lefschetz trace formula for open adic spaces
over an algebraically closed non-archimedean field.
First recall the Lefschetz-Verdier trace formula for schemes.
Let $X$ be a scheme which is separated smooth of finite type over an algebraically closed field $k$,
$X\hooklongrightarrow \overline{X}$ a dense compactification
and $\overline{f}\colon \overline{X}\longrightarrow \overline{X}$ a $k$-morphism 
which induces a proper $k$-morphism $f\colon X\longrightarrow X$.
Assume for simplicity that the fixed scheme $\Fix f$ defined by the cartesian diagram
\[
 \xymatrix{%
 \Fix f\ar[r]\ar[d]&X\ar[d]^-{\text{diagonal}}\\
 X\ar[r]^-{f\times\id}&X\times X
 }
\]
is discrete. Then, for a prime $\ell$ which is invertible in $k$,
the alternating sum of the traces $\sum_{i}(-1)^i\Tr(f^*;H^i_c(X,\Q_\ell))$ is equal to
$\#\Fix f+\sum_{D\in \pi_0(\Fix \overline{f}\cap (\overline{X}\setminus X))}\loc_D(f)$,
where $\#\Fix f$ is the number of fixed points by $f$ counted with multiplicity, $\pi_0(-)$ denotes
the set of connected components and $\loc_D(f)$ denotes the ``local term at $D$''
(\cf \cite[\S 1.2]{MR1431137}).
In particular, if $\Fix \overline{f}\subset X$ then $\sum_{i}(-1)^i\Tr(f^*;H^i_c(X,\Q_\ell))$ coincides 
with $\#\Fix f$, which gives a nice fixed point formula.

It seems natural to expect the similar formula for  adic spaces.
However, there is an obvious counterexample. Let $k$ be an algebraically closed non-archimedean field and
$\mathbb{D}^1=\Spa(k\langle T\rangle,k\langle T\rangle^\circ)$ the unit disk.
We can compactify it by taking its closure $\overline{\mathbb{D}^1}$ in $(\A^1)^\ad$.
Consider the isomorphism $\A^1\longrightarrow \A^1$ given by $T\longmapsto T+1$, which induces
the isomorphisms $f\colon \mathbb{D}^1\longrightarrow \mathbb{D}^1$ and 
$\overline{f}\colon \overline{\mathbb{D}^1}\longrightarrow \overline{\mathbb{D}^1}$.
Since $\Fix f=\Fix \overline{f}=\varnothing$, we expect to have
$\sum_{i}(-1)^i\Tr(f^*;H^i_c(\mathbb{D}^1,\Q_\ell))=0$.
Nevertheless the left hand side is equal to $1$, and thus the analogue of the Lefschetz-Verdier trace formula do not
hold.

Actually, this phenomenon has already been observed by Fujiwara \cite{MR1431137} and Huber \cite{MR1856262}.
Fujiwara proved his topological Lefschetz trace formula under the condition that there exists no topological fixed point
on the boundary. 
Huber established a trace formula for open curves, which says that
the alternating sum of the traces on the cohomology of an open adic curve $X$ is the sum of the number of fixed points
and the contribution at each set-theoretical fixed point on the boundary $X^c\setminus X$, where
$X^c$ denotes the universal compactification of $X$.
Our main theorem is also in this line. Here we will give a slightly simplified statement.
As above, let $k$ be an algebraically closed non-archimedean field and $\ell$ a prime which is invertible in the
residue field of $k$.
Let $X$ be a purely $d$-dimensional adic space which is separated smooth of finite type over $k$,
$X\hooklongrightarrow \overline{X}$ a dense compactification
and $\overline{f}\colon \overline{X}\longrightarrow \overline{X}$ a $k$-morphism 
which induces a proper $k$-morphism $f\colon X\longrightarrow X$.

\begin{thm}\label{thm:main-thm}
 Assume that for every $x\in \overline{X}\setminus X$, the points $x$ and $\overline{f}(x)$
 are distinct.
 Then we have
 \[
  \Tr\bigl(f^*;R\Gamma_c(X,\Z/\ell^n\Z)\bigr)=\#\Fix f.
 \]
 If moreover the characteristic of $k$ is $0$, then 
 \[
  \sum_{i=0}^{2d}(-1)^i\Tr\bigl(f^*;H^i_c(X,\Q_\ell)\bigr)=\#\Fix f.
 \]
\end{thm}

Since we have no intersection theory for adic spaces yet (at least the author do not know),
we need to clarify the meaning of ``the number of fixed points'' $\#\Fix f$. 
The definition is given by using cohomology theory (see Definition \ref{defn:Fix}).

The statement above is similar to \cite[Theorem 2.2.8]{MR1431137}, but our theorem is valid for a non-algebraizable case.
Our proof is very different from that in \cite{MR1431137}; we use neither formal geometry nor the Lefschetz-Verdier
trace formula for schemes. Our proof is purely rigid-geometric.
The first step of our proof is to observe that why the proof of the Lefschetz-Verdier
trace formula in \cite[Expos\'e III]{SGA5} cannot be applied to the case of adic spaces;
actually, the only obstruction is the failure of the K\"unneth formula for push-forward 
$(Rf_*F)\Lboxtimes (Rg_*G)\cong R(f\times g)_*(F\Lboxtimes G)$
(\cf Remark \ref{rem:isom-failure}). Therefore, our main strategy is to find a suitable isomorphism
induced by the K\"unneth homomorphism $(Rf_*F)\Lboxtimes (Rg_*G)\longrightarrow R(f\times g)_*(F\Lboxtimes G)$
by using the assumption in Theorem \ref{thm:main-thm},
so that the analogous proof as in \cite[Expos\'e III]{SGA5} works.
This idea is also useful for finding other trace formulas than Theorem \ref{thm:main-thm}. For example,
we will give another formula for open curves which is very similar to the formula by Huber,
and a partial generalization of Fujiwara's trace formula for contracting correspondences
to the non-algebraizable case.

The author's main motivation for this work is to establish the Lefschetz trace formula which is 
applicable to the Rapoport-Zink towers. For a classical case, there are such works by Faltings \cite{MR1302321}
and Strauch \cite{MR2383890}; the former is on the Drinfeld tower and the latter is on the Lubin-Tate tower.
As an application of Theorem \ref{thm:main-thm}, we will establish the Lefschetz trace formula for formal schemes
(Theorem \ref{thm:LTF-formal}),
which is applicable to the Rapoport-Zink tower for $\mathrm{GSp}(4)$ considered in \cite{RZ-GSp4}.
In a forthcoming paper, the author will use this formula to observe that
the local Jacquet-Langlands correspondence for $\mathrm{GSp}(4)$ appears in the $\ell$-adic cohomology
of the Rapoport-Zink tower mentioned above. 
He also hopes that there are a few more Rapoport-Zink towers
to which our trace formula can apply.

We sketch the outline of the paper. In Section 2, we will consider the action $\gamma^*$ of a correspondence $\gamma$
on the \'etale cohomology of adic spaces, and give the definition of ``the number of fixed points'' $\#\Fix \gamma$.
We also show that it is \'etale local and compatible with the comparison functor.
These properties justify our definition, though it is cohomological and far from geometric.
In Section 3, we will prove our main theorem. First we discuss the K\"unneth formula for $R\Gamma$,
which leads us to a weaker form of the Lefschetz trace formula (Proposition \ref{prop:LTF-unlocalized}).
Next we refine this weaker version by using our assumption on points of the boundary to get our main theorem.
We also remark on the trace formula for open curves.
In Section 4, we prove the Lefschetz trace formula for formal schemes and give some interesting examples.
In Section 5, we give a simple trace formula for a morphism which is contracting near fixed points.
It is a partial generalization of a result of Fujiwara \cite[Theorem 3.2.4]{MR1431137}.

\bigbreak

\noindent{\bfseries Acknowledgment}\quad
The author would like to thank an anonymous referee for valuable comments.
This work was supported by JSPS KAKENHI Grant Number 21740022.

\bigbreak

\noindent{\bfseries Notation}\quad
Let $k$ be an algebraically closed non-archimedean field (\cf \cite[Definition 1.1.3]{MR1734903})
and denote its valuation ring by $k^+$. 
Put $S=\Spa(k,k^+)$.
Fix a prime $\ell$ which is invertible in $k^+$ and put $\Lambda=\Z/\ell^n\Z$
for an integer $n\ge 1$.

Every sheaf and cohomology are considered in the \'etale topology.
We simply write $f^!$ for the functor $R^+f^!$ introduced in \cite[Theorem 7.1.1]{MR1734903}.

\section{Correspondences on adic spaces}
\subsection{Trace maps and Gysin maps}
 Let $X$ be a purely $d$-dimensional adic space which is separated, locally of finite type and taut over $S$.
 Assume that $X$ is generically smooth over $S$, namely, the dimension of the singular locus $Z$ of $X$ is 
 strictly less than $d$.
 We will construct the trace map $\Tr_X\colon H^{2d}_c(X,\Lambda(d))\longrightarrow \Lambda$ for such $X$.

 Note that the complement $U$ of $Z$ is taut \cite[Lemma 5.1.4 i)]{MR1734903},
 thus we may define the compactly supported cohomology of $X$, $Z$ and $U$.
 By the exact sequence
 \[
  H^{2d-1}_c\bigl(Z,\Lambda(d)\bigr)\longrightarrow H^{2d}_c\bigl(U,\Lambda(d)\bigr)\longrightarrow 
  H^{2d}_c\bigl(X,\Lambda(d)\bigr)\longrightarrow H^{2d}_c\bigl(Z,\Lambda(d)\bigr)
 \]
 and the vanishing $H^{2d-1}_c(Z,\Lambda(d))=H^{2d}_c(Z,\Lambda(d))=0$
 (\cite[Proposition 5.5.8, Corollary 1.8.8]{MR1734903}),
 the canonical homomorphism $H^{2d}_c(U,\Lambda(d))\longrightarrow H^{2d}_c(X,\Lambda(d))$ is an isomorphism.
 
 On the other hand, since $U$ is smooth over $S$,
 the trace map $\Tr_U\colon H^{2d}_c(U,\Lambda(d))\longrightarrow \Lambda$ has already been constructed in
 \cite[Theorem 7.3.4]{MR1734903}.
 We will define $\Tr_X$ as the composite
 \[
  H^{2d}_c\bigl(X,\Lambda(d)\bigr)\yleftarrow{\cong} H^{2d}_c\bigl(U,\Lambda(d)\bigr)\yrightarrow{\Tr_U}\Lambda.
 \]
 Since $H^i_c(X,\Lambda(d))=0$ for $i>2d$, $\Tr_X$ induces the map
 $R\Gamma_c(X,\Lambda(d)[2d])\longrightarrow \Lambda$,
 which is also denoted by $\Tr_X$.

\begin{prop}\label{prop:trace-etale}
 Let $X$, $X'$ be purely $d$-dimensional adic spaces which are separated, locally of finite type, generically smooth
 and taut over $S$.
 Let $\pi\colon X'\longrightarrow X$ be an \'etale $S$-morphism between them. Then the composite
 $H^{2d}_c(X',\Lambda(d))\yrightarrow{\pi_*} H^{2d}_c(X,\Lambda(d))\yrightarrow{\Tr_X}\Lambda$ is equal to $\Tr_{X'}$.
\end{prop}

\begin{prf}
 First assume that $X$ is smooth over $S$. Then the claim follows from the characterizing properties of
 the trace morphisms ((Var 3) and (Var 4) in \cite[Theorem 7.3.4]{MR1734903}).
 In the general case, let $Z$ (resp.\ $Z'$) be the singular locus of $X$ (resp.\ $X'$) and put 
 $U=X\setminus Z$ (resp.\ $U'=X'\setminus Z'$).
 Note that we have $\pi^{-1}(U)\subset U'$, for $\pi$ is \'etale. 
 Thus we have the following commutative diagram:
 \[
  \xymatrix{%
 H^{2d}_c\bigl(X',\Lambda(d)\bigr)\ar@{=}[d]& H^{2d}_c\bigl(U',\Lambda(d)\bigr)\ar[l]_-{\cong}\ar[rr]^-{\Tr_{U'}}&& \Lambda\ar@{=}[d]\\
 H^{2d}_c\bigl(X',\Lambda(d)\bigr)\ar[d]^-{\pi_*}&H^{2d}_c\bigl(\pi^{-1}(U),\Lambda(d)\bigr)\ar[d]^-{\pi_*}\ar[l]_-{(*)}\ar[rr]^-{\Tr_{\pi^{-1}(U)}}\ar[u]
 &&\Lambda\ar@{=}[d]\\
 H^{2d}_c\bigl(X,\Lambda(d)\bigr)&H^{2d}_c\bigl(U,\Lambda(d)\bigr)\ar[l]_-{\cong}\ar[rr]_-{\Tr_U}&& \Lambda\lefteqn{.}
 }
 \]
 Here the homomorphism $(*)$ is an isomorphism, since the closed adic subspace $\pi^{-1}(Z)$ of $X'$,
 which is \'etale over $Z$, has dimension less than $d$.
 Therefore the claim immediately follows from the diagram above.
\end{prf}

\begin{prop}\label{prop:trace-GAGA}
 Let $X$ be a purely $d$-dimensional scheme which is separated of finite type over $k$ and
 $X^\ad=X\times_{\Spec k}S$ the associated adic space.
 Assume that $X$ is generically smooth over $k$.
 Then the composite of the canonical comparison map $H^{2d}_c(X,\Lambda(d))\longrightarrow H^{2d}_c(X^\ad,\Lambda(d))$ and $\Tr_{X^\ad}\colon H^{2d}_c(X^\ad,\Lambda(d))\longrightarrow \Lambda$ is equal to the trace map
 $\Tr_X\colon H^{2d}_c(X,\Lambda(d))\longrightarrow \Lambda$ for $X$, whose definition is similar to that of
 $\Tr_{X^\ad}$.
\end{prop}

First we consider the smooth case.

\begin{prop}\label{prop:trace-GAGA-smooth}
 Let $f\colon X\longrightarrow Y$ be a separated smooth morphism of finite type with relative dimension $d$
 between $k$-schemes of finite type.
 Let $f^\ad\colon X^\ad\longrightarrow Y^\ad$ be the induced morphism of adic spaces.
 Denote the natural morphisms of sites $(X^\ad)_\et\longrightarrow X_\et$ and 
 $(Y^\ad)_\et\longrightarrow Y_\et$ by $\varepsilon$.
 Then the following diagram is commutative:
 \[
 \xymatrix{%
 \varepsilon^*Rf_!\Lambda(d)[2d]\ar[rr]^-{\varepsilon^*\Tr_f}\ar[d]&& \Lambda\ar@{=}[d]\\
 Rf^\ad_!\Lambda(d)[2d]\ar[rr]^-{\varepsilon^*\Tr_{f^\ad}}&& \Lambda\lefteqn{.}
 }
 \]
 The left vertical arrow is defined in \cite[Theorem 5.7.2]{MR1734903}.
\end{prop}

\begin{prf}
 By using \cite[Expos\'e XVIII, Lemme 2.2]{SGA4} and \cite[Lemma 7.3.5]{MR1734903}, we have only to
 consider the case where $X=\A^m_Y$ or the case where $f$ is \'etale.
 The former case immediately follows from the construction of $\Tr_{f^\ad}$
 (\cite[proof of Theorem 7.3.4]{MR1734903}). For the latter case, every morphism is defined by the 
 adjointness and the commutativity is formal.
\end{prf}

\begin{prf}[of Proposition \ref{prop:trace-GAGA}]
 Let $Z$ be the singular locus of $X$ and put $U=X\setminus Z$. Then, by Proposition \ref{prop:trace-GAGA-smooth},
 the following diagram is commutative:
 \[
  \xymatrix{%
 H^{2d}_c\bigl(X,\Lambda(d)\bigr)\ar[d]& H^{2d}_c\bigl(U,\Lambda(d)\bigr)\ar[d]\ar[l]_-{\cong}\ar[rr]^-{\Tr_U}&& \Lambda\ar@{=}[d]\\
 H^{2d}_c\bigl(X^\ad,\Lambda(d)\bigr)& H^{2d}_c\bigl(U^\ad,\Lambda(d)\bigr)\ar[l]_-{(*)}\ar[rr]^-{\Tr_{U^\ad}}&& \Lambda\lefteqn{.} 
 }
 \]
 Since the dimension of $Z^\ad$ is less than $d$, the map $(*)$ is an isomorphism.
 Moreover, in the same way as Proposition \ref{prop:trace-etale}, we can prove that
 the composite of 
 $H^{2d}_c(X^\ad,\Lambda(d))\yleftarrow{\cong}H^{2d}_c(U^\ad,\Lambda(d))\yrightarrow{\Tr_{U^\ad}}\Lambda$
 coincides with $\Tr_{X^\ad}$. Thus we have the desired compatibility.
\end{prf}

Let $X$ (resp.\ $Y$) be a purely $d$-dimensional (resp.\ $d'$-dimensional) adic space which is separated,
locally of finite type and taut over $S$. Put $c=d-d'$. Assume $X$ (resp.\ $Y$) is smooth (resp.\ generically smooth)
over $S$. 
Let $f\colon Y\longrightarrow X$ be an $S$-morphism between them.
Let us denote the structure map of $X$ (resp.\ $Y$) by $a\colon X\longrightarrow S$
(resp.\ $b\colon Y\longrightarrow S$). By the construction above, we have 
$\Tr_X\colon Ra_!\Lambda(d)[2d]\longrightarrow \Lambda$ and $\Tr_Y\colon Rb_!\Lambda(d')[2d']\longrightarrow \Lambda$.
By the adjointness, these correspond to the maps $\Gys_a\colon \Lambda\longrightarrow a^!\Lambda(-d)[-2d]$ and 
$\Gys_b\colon \Lambda\longrightarrow b^!\Lambda(-d')[-2d']$. Since $a$ is smooth, $\Gys_a$ is an isomorphism 
(\cite[Theorem 7.5.3]{MR1734903}).
Therefore we have an isomorphism $b^!\Lambda(-d')[-2d']=f^!a^!\Lambda(-d')[-2d']\cong f^!\Lambda(c)[2c]$, and
finally we obtain a map $\Gys_f\colon \Lambda\longrightarrow f^!\Lambda(c)[2c]$, which is called
the Gysin map associated with $f$. Since $\Hom(\Lambda,f^!\Lambda(c)[2c])=H^{2c}(Y,f^!\Lambda(c))$,
it gives an element of $H^{2c}(Y,f^!\Lambda(c))$.

If moreover $f$ is proper, then $f$ is naturally decomposed as
$Y\yrightarrow{f'} f(Y)\yrightarrow{i} X$, where $f'$ is proper and $i$ is a closed immersion
(here $f(Y):=(X,f(Y))$ is a pseudo-adic space; \cf \cite[\S1.10]{MR1734903}). Thus we have the map
\begin{align*}
 H^{2c}\bigl(Y,f^!\Lambda(c)\bigr)&=H^{2c}\bigl(f(Y),Rf'_!f'^!i^!\Lambda(c)\bigr)\\
 &\yrightarrow{\adj}H^{2c}\bigl(f(Y),i^!\Lambda(c)\bigr)=H^{2c}_{f(Y)}\bigl(X,\Lambda(c)\bigr).
\end{align*}
We denote the image of $\Gys_f$ under this map by $\cl(f)$,
and call it the cohomology class associated with $f$.

\begin{prop}\label{prop:Gysin-etale}
 In the setting above (we do not assume that $f$ is proper), consider the following commutative diagram:
 \[
 \xymatrix{%
 Y'\ar[r]^-{g}\ar[d]^-{\pi'}&X'\ar[d]^-{\pi}\\
 Y\ar[r]^-{f}&X\lefteqn{,}
 }
 \]
 where $\pi$ and $\pi'$ are \'etale. 
 Then the image of $\Gys_f$ under the map
 \[
  \pi^*\colon H^{2c}\bigl(Y,f^!\Lambda(c)\bigr)\yrightarrow{\pi'^*} H^{2c}\bigl(Y',\pi'^*f^!\Lambda(c)\bigr)
 =H^{2c}\bigl(Y',g^!\Lambda(c)\bigr)
 \]
 coincides with $\Gys_g$.  
 If moreover the diagram above is cartesian and $f$ is proper, then the image of $\cl(f)$ under the map
 $\pi^*\colon H^{2c}_{f(Y)}(X,\Lambda(c))\longrightarrow H^{2c}_{g(Y')}(X',\Lambda(c))$ coincides with $\cl(g)$.
\end{prop}

\begin{prf}
 First we will prove $\pi^*\Gys_f=\Gys_g$.
 Denote the structure map of $X'$ (resp.\ $Y'$) by $a'$ (resp.\ $b'$).
 By Proposition \ref{prop:trace-etale} and the adjointness, we have the following
 commutative diagrams:
 \[
  \xymatrix{%
 \Lambda\ar[rr]^-{\Gys_{a'}}_-{\cong}\ar@{=}[d]&& a'^!\Lambda(-d)[-2d]\ar@{=}[d]\\
 \pi^*\Lambda\ar[rr]^-{\pi^*\Gys_a}_-{\cong}&& \pi^*a^!\Lambda(-d)[-2d]\lefteqn{,}
 }\qquad
  \xymatrix{%
 \Lambda\ar[rr]^-{\Gys_{b'}}\ar@{=}[d]&& b'^!\Lambda(-d')[-2d']\ar@{=}[d]\\
 \pi'^*\Lambda\ar[rr]^-{\pi'^*\Gys_b}&& \pi'^*b^!\Lambda(-d')[-2d']\lefteqn{.}
 }
 \]
 Namely, we have $\Gys_{a'}=\pi^*\Gys_a$ and $\Gys_{b'}=\pi'^*\Gys_b$.
 On the other hand, by the definition, $\Gys_f=f^!\Gys_a^{-1}(c)[2c]\circ \Gys_b$ and 
 $\Gys_g=g^!\Gys_{a'}^{-1}(c)[2c]\circ \Gys_{b'}$. Therefore we have
 \begin{align*}
 \pi^*\Gys_f&=\pi^*f^!\Gys_a^{-1}(c)[2c]\circ \pi^*\Gys_b=g^!\pi'^*\Gys_a^{-1}(c)[2c]\circ \pi^*\Gys_b\\
  &=g^!\Gys_{a'}^{-1}(c)[2c]\circ \Gys_{b'}=\Gys_g,
 \end{align*}
 as desired.

 Assume that the diagram in the proposition is cartesian and $f$ is proper.
 To prove $\pi^*\cl(f)=\cl(g)$, it suffices to observe the commutativity of the diagram below:
 \[
 \xymatrix{%
 H^{2c}\bigl(Y,f^!\Lambda(c)\bigr)\ar[r]\ar[d]^-{\pi'^*}&H^{2c}_{f(Y)}\bigl(X,\Lambda(c)\bigr)\ar[d]^-{\pi^*}\\
 H^{2c}\bigl(Y',g^!\Lambda(c)\bigr)\ar[r]&H^{2c}_{g(Y')}\bigl(X',\Lambda(c)\bigr)\lefteqn{.}
 }
 \]
 Let $Y'\yrightarrow{g'} g(Y')\yrightarrow{i'} X$ be the factorization of $g$.
 Note that the following diagram is cartesian due to \cite[Lemma 3.9 (i)]{MR1306024}:
 \[
  \xymatrix{%
 Y'\ar[d]^-{\pi'}\ar[r]^-{g'}& g(Y')\ar[r]^-{i'}\ar[d]^-{\pi}& X'\ar[d]^-{\pi}\\
 Y\ar[r]^-{f'}& f(Y)\ar[r]^-{i}& X\lefteqn{.}
 }
 \]
 For every object $L$ of $D^+(f(Y),\Lambda)$, it is straightforward to check the commutativity of 
 the diagram below:
 \[
 \xymatrix{%
 Rf'_!f'^!L\ar[r]^-{\adj}\ar[d]^-{\adj}& L\ar[d]^-{\adj}\\
 R\pi_*\pi^*Rf'_!f'^!L\ar[r]^-{\adj}\ar@{<->}[d]_-{\cong}^-{\mathrm{base\ change}}&
 R\pi_*\pi^*L\ar@{=}[dd]\\
 R\pi_*Rg'_!\pi'^*f'^!L\ar@{=}[d]\\
 R\pi_*Rg'_!g'^!\pi^*L\ar[r]^-{\adj}& R\pi_*\pi^*L\lefteqn{.}
 }
 \]
 Setting $L=i^!\Lambda(c)$ and taking $H^{2c}(f(Y),-)$, we obtain the desired commutativity.
\end{prf}

Next we will prove a comparison result for the Gysin maps and the cohomology classes
associated with proper morphisms.
Let $X$ (resp.\ $Y$) be a purely $d$-dimensional (resp.\ $d'$-dimensional) scheme which is separated
of finite type over $k$. Assume $X$ (resp.\ $Y$) is smooth (resp.\ generically smooth) over $k$ and put $c=d-d'$. 
Let $f\colon Y\longrightarrow X$ be a morphism of finite type over $k$. Then, by the same way as above,
we may define the Gysin map $\Gys_f\colon \Lambda\longrightarrow f^!\Lambda(c)[2c]$.
If moreover $f$ is proper, we can also define the cohomology class $\cl(f)\in H^{2c}_{f(Y)}(X,\Lambda(c))$
associated with $f$.
Let $f^{\ad}\colon Y^{\ad}\longrightarrow X^{\ad}$ be the morphism of adic spaces associated with $f$
and denote the natural morphisms of locally ringed spaces $X^\ad\longrightarrow X$ and $Y^\ad\longrightarrow Y$
by $\varepsilon$.
Note that $X^{\ad}$, $Y^\ad$ and $f^\ad$ satisfy all the assumptions above,
thus we may define the Gysin map $\Gys_{f^\ad}$,
and the class $\cl(f^{\ad})$ if $f$ is proper.

\begin{prop}\label{prop:Gysin-GAGA}
 The image of $\Gys_f$ under the map
 \[
  H^{2c}\bigl(Y,f^!\Lambda(c)\bigr)\yrightarrow{\varepsilon^*} H^{2c}\bigl(Y,\varepsilon^*f^!\Lambda(c)\bigr)
 \yrightarrow[(*)]{\cong} H^{2c}\bigl(Y,f^{\ad!}\varepsilon^*\Lambda(c)\bigr)=H^{2c}\bigl(Y,f^{\ad!}\Lambda(c)\bigr)
 \]
 coincides with $\Gys_{f^\ad}$ (for the construction of the isomorphism $(*)$,
 see \cite[Proposition 4.37]{formalnearby}). 
 If moreover $f$ is proper, the image of $\cl(f)$ under the canonical map 
 $\varepsilon^*\colon H^{2d}_{f(Y)}(X,\Lambda(d))\longrightarrow H^{2d}_{f^{\ad}(Y^{\ad})}(X^{\ad},\Lambda(d))$
 coincides with $\cl(f^{\ad})$
 (note that $f^\ad(Y^\ad)=\varepsilon^{-1}(f(Y))$ due to \cite[Lemma 3.9 (ii)]{MR1306024}).
\end{prop}

\begin{prf}
 First we will observe $\varepsilon^*\Gys_f=\Gys_{f^\ad}$.
 Let us denote the structure morphism of $X$ (resp.\ $Y$) by $a$ (resp.\ $b$).
 By Proposition \ref{prop:trace-GAGA} and the adjointness,
 we have the commutativity of the left diagram below. On the other hand, by the definition of $\Gys_b$,
 the right diagram below is also commutative.
 \[
  \xymatrix{%
 \varepsilon^*\Lambda\ar[r]^-{\Gys_{b^\ad}}\ar[d]^-{\adj}& b^{\ad!}\varepsilon^*\Lambda(-d')[-2d']\ar@{=}[dd]\\
 b^{\ad !}Rb_!^\ad\varepsilon^*\Lambda\ar@{<->}[d]_-{\cong}^-{\text{b.c.}}\\
 b^{\ad !}\varepsilon^*Rb_!\Lambda\ar[r]^-{b^{\ad !}\!\varepsilon^*\!\Tr_b}& b^{\ad!}\varepsilon^*\Lambda(-d')[-2d']
 \lefteqn{,}}
 \quad
 \xymatrix{%
 \varepsilon^*\Lambda\ar[r]^-{\varepsilon^*\Gys_b}\ar[d]^-{\adj}& \varepsilon^*b^!\Lambda(-d')[-2d']\ar@{=}[d]\\
 \varepsilon^*b^!Rb_!\Lambda\ar[r]^-{\varepsilon^*\!b^!\!\Tr_b}\ar[d]_-{\cong}^-{\text{b.c.}}& \varepsilon^*b^!\Lambda(-d')[-2d']\ar[d]_-{\cong}^-{\text{b.c.}}\\
 b^{\ad !}\varepsilon^*Rb_!\Lambda\ar[r]^-{b^{\ad !}\!\varepsilon^*\!\Tr_b}& b^{\ad!}\varepsilon^*\Lambda(-d')[-2d']
 \lefteqn{.}}
 \]
 Moreover, it is easy to show that the following diagram is commutative:
 \[
  \xymatrix{%
 \varepsilon^*\Lambda\ar[r]^-{\adj}\ar[d]^-{\adj}& b^{\ad !}Rb^\ad_!\varepsilon^*\Lambda\\
 \varepsilon^*b^!Rb_!\Lambda\ar[r]^-{\text{b.c.}}& b^{\ad !}\varepsilon^*Rb_!\Lambda\ar[u]_-{\text{b.c.}}\lefteqn{.}
 }
 \]
 By these three commutative diagrams, we have the following commutative diagram:
 \[
  \xymatrix{%
 \varepsilon^*\Lambda\ar[rr]^-{\Gys_{b^\ad}}\ar@{=}[d]&& \varepsilon^*b^!\Lambda\ar[d]^-{\text{b.c.}}\\
 \varepsilon^*\Lambda\ar[rr]^-{\varepsilon^*\Gys_b}&& b^{\ad !}\varepsilon^*\Lambda\lefteqn{.}
 }
 \]
 Namely, we have $\Gys_{b^\ad}=\varepsilon^*\Gys_b$. Similarly we have $\Gys_{a^\ad}=\varepsilon^*\Gys_a$.
 Therefore, in the same way as in the proof of Proposition \ref{prop:Gysin-etale}, we can obtain
 $\varepsilon^*\Gys_f=\Gys_{f^\ad}$.

 Assume that $f$ is proper. To prove $\varepsilon^*\cl(f)=\cl(f^\ad)$, it suffices to show
 the commutativity of the diagram below:
 \[
 \xymatrix{%
 H^{2c}\bigl(Y,f^!\Lambda(c)\bigr)\ar[r]\ar[d]^-{\varepsilon^*}&H^{2c}_{f(Y)}\bigl(X,\Lambda(c)\bigr)\ar[d]^-{\varepsilon^*}\\
 H^{2c}\bigl(Y^\ad,f^{\ad!}\Lambda(c)\bigr)\ar[r]&H^{2c}_{f^\ad(Y^\ad)}\bigl(X^\ad,\Lambda(c)\bigr)\lefteqn{.}
 }
 \]
 Put $Z=f(Y)$ and endow it with the structure of a reduced closed subscheme of $X$.
 Then $f$ factors as $Y\yrightarrow {f'} Z\stackrel{i}{\hooklongrightarrow} X$. As we mentioned above,
 a closed subset $f^\ad(Y^\ad)$ of $X^\ad$ coincides with $Z^\ad$. Therefore, by \cite[Corollary 2.3.8]{MR1734903},
 $i^\ad\colon Z^\ad \hooklongrightarrow X^\ad$ induces an equivalence between the \'etale topos of $Z^\ad$ and
 that of $f^\ad(Y^\ad)$. Hence the diagram above can be identified with the following diagram:
 \[
 \xymatrix{%
 H^{2c}\bigl(Y,f^!\Lambda(c)\bigr)\ar[r]\ar[d]^-{\varepsilon^*}&H^{2c}\bigl(Z,i^!\Lambda(c)\bigr)\ar[d]^-{\varepsilon^*}\\
 H^{2c}\bigl(Y^\ad,f^{\ad!}\Lambda(c)\bigr)\ar[r]&H^{2c}\bigl(Z^\ad,i^{\ad!}\Lambda(c)\bigr)\lefteqn{.}
 }
 \]
 Now we can show the commutativity in the same way as in the proof of Proposition \ref{prop:Gysin-etale}.
 This completes the proof.
\end{prf}

\subsection{Correspondences}\label{subsec:corr}
Let $X$ and $\Gamma$ be purely $d$-dimensional adic spaces which are separated, locally of finite type and taut
over $S$.
Assume that $X$ (resp.\ $\Gamma$) is smooth (resp.\ generically smooth) over $S$.
Let $\gamma\colon \Gamma\longrightarrow X\times_SX$ be a morphism over $S$ and 
put $\gamma_i=\pr_i\circ \gamma$. 

Then we may apply the construction in the previous subsection and have the cohomology class
$\Gys_\gamma\in H^{2d}(\Gamma,\gamma^!\Lambda(d))$.
Using it, we will define ``the number of points fixed by $\gamma$''.

\begin{defn}\label{defn:Fix}
 Consider the following cartesian diagram in which $\delta\colon X\longrightarrow X\times_SX$ denotes
 the diagonal morphism:
 \[
  \xymatrix{%
 \Fix \gamma\ar[r]^-{\gamma_0}\ar[d]& X\ar[d]^-{\delta}\\
 \Gamma\ar[r]^-{\gamma}& X\times_SX\lefteqn{.}
 }
 \]
 Let $D$ be an open and closed subset of $\Fix \gamma$ which is proper over $S$,
 and denote the open and closed immersion $D\hooklongrightarrow \Fix \gamma$ by $j$.
 Put $j_0=\gamma_0\circ j$.
 Then we have the canonical maps 
 \begin{align*}
  H^{2d}\bigl(\Gamma,\gamma^!\Lambda(d)\bigr)&\yrightarrow{\delta^*}
  H^{2d}\bigl(\Fix \gamma,\gamma_0^!\Lambda(d)\bigr)\longrightarrow
  H^{2d}\bigl(D,j^*\gamma_0^!\Lambda(d)\bigr)
  =H^{2d}\bigl(D,j_0^!\Lambda(d)\bigr)\\
  &=H^{2d}_c\bigl(X,Rj_{0!}j_0^!\Lambda(d)\bigr)
  \yrightarrow{\adj}H^{2d}_c\bigl(X,\Lambda(d)\bigr)\yrightarrow{\Tr}\Lambda.
 \end{align*}
 We denote the image of $\Gys_\gamma$ under these maps by $\#\Fix_D\gamma$.
 If $D=\Fix\gamma$, we write $\#\Fix\gamma$ for $\#\Fix_D\gamma$.
\end{defn}

\begin{rem}\label{rem:Fix-via-cl}
 Assume that $\gamma$ is proper. Then $\#\Fix\gamma$ can be calculated from $\cl(\gamma)$ as follows.
 Denote the inverse image of $\gamma(\Gamma)$ under $\delta$ by $\Delta_X\cap \gamma(\Gamma)$.
 Since we are implicitly assuming that $\Fix\gamma$ is proper over $S$, the pseudo-adic space
 $(X,\Delta_X\cap \gamma(\Gamma))$ is proper over $S$. Thus we have the following natural maps:
 \[
  H^{2d}_{\gamma(\Gamma)}\bigl(X\times_SX,\Lambda(d)\bigr)\yrightarrow{\delta^*}
 H^{2d}_{\Delta_X\cap \gamma(\Gamma)}\bigl(X,\Lambda(d)\bigr)\longrightarrow
 H^{2d}_c\bigl(X,\Lambda(d)\bigr)\yrightarrow{\Tr_X}\Lambda.
 \]
 The image of $\cl(\gamma)$ under these maps coincides with $\#\Fix \gamma$.
\end{rem}

\begin{rem}\label{rem:Fix-l-adic}
 Obviously the construction of $\#\Fix \gamma$ is compatible with a change of $\Lambda$.
 Namely, if we denote $\#\Fix \gamma$ for $\Lambda=\Z/\ell^n\Z$ by $\#_{\ell^n}\!\Fix \gamma$,
 then $(\#_{\ell^n}\!\Fix \gamma)_{n\ge 1}$ gives an element of $\Z_\ell=\varprojlim_{n}\Z/\ell^n\Z$.
 We also denote it by $\#\Fix \gamma$.
\end{rem}

\begin{exa}\label{exa:self-mor}
 Let $f\colon X\longrightarrow X$ be a morphism over $S$.
 Then $\gamma_f=f\times \id\colon X\longrightarrow X\times_SX$ lies in the situation above.
 In this case, $\Fix \gamma_f$ is a closed adic subspace of $X$.
 We simply write $\Fix f$, $\#\Fix_Df$ and $\#\Fix f$
 for $\Fix \gamma_f$, $\#\Fix_D\gamma_f$ and $\#\Fix\gamma_f$, respectively.
 Denote the image of $\gamma_f$ by $\Gamma_f$; note that it has the natural structure of a closed adic subspace
 of $X\times_SX$ since $\gamma_f$ is a closed immersion.
\end{exa}

By the results in the previous subsection, we can prove that the number $\#\Fix_D\gamma$ is 
\'etale local and compatible with the comparison functor:

\begin{prop}\label{prop:Fix-etale}
 Let $\gamma\colon \Gamma\longrightarrow X\times_SX$ be as above and
 $\gamma'\colon \Gamma'\longrightarrow X'\times_SX'$ be another morphism
 satisfying the conditions above.
 Let $\pi\colon \Gamma'\longrightarrow \Gamma$ and $\pi'\colon X'\longrightarrow X$ be
 \'etale morphisms over $S$ such that $\gamma\circ \pi=(\pi'\times \pi')\circ \gamma'$.
 Let $D$ (resp.\ $D'$) be an open and closed subset of $\Fix \gamma$
 (resp.\ $\Fix\gamma'$) which is proper over $S$.
 Assume that $\pi$ induces an isomorphism from $D'$ to $D$.
 Then we have $\#\Fix_D\gamma=\#\Fix_{D'}\gamma'$.
\end{prop}

\begin{prf}
 We have the following commutative diagram:
 \[
  \xymatrix{%
 H^{2d}\bigl(\Gamma,\gamma^!\Lambda(d)\bigr)\ar[r]^-{\delta^*}\ar[d]^-{\pi^*}&
 H^{2d}\bigl(\Fix\gamma,\gamma_0^!\Lambda(d)\bigr)\ar[r]\ar[d]^-{\pi^*}&
 H^{2d}\bigl(D,j_0^!\Lambda(d)\bigr)\ar[d]^-{\pi^*}\\
 H^{2d}\bigl(\Gamma',\gamma'^!\Lambda(d)\bigr)\ar[r]^-{\delta'^*}&
 H^{2d}\bigl(\Fix\gamma',\gamma_0'^!\Lambda(d)\bigr)\ar[r]&
 H^{2d}\bigl(D',j_0'^!\Lambda(d)\bigr)\lefteqn{.}
 }
 \]
 Thus, by Proposition \ref{prop:trace-etale} and Proposition \ref{prop:Gysin-etale},
 it suffices to show the commutativity of the diagram below:
 \[
  \xymatrix{%
 H^{2d}\bigl(D,j_0^!\Lambda(d)\bigr)\ar@{=}[r]\ar[d]^-{\pi^*}&
 H^{2d}_c\bigl(X,Rj_{0!}j_0^!\Lambda(d)\bigr)\ar[r]& H^{2d}_c\bigl(X,\Lambda(d)\bigr)\\
 H^{2d}\bigl(D',j_0'^!\Lambda(d)\bigr)\ar@{=}[r]&
 H^{2d}_c\bigl(X',Rj'_{0!}j_0'^!\Lambda(d)\bigr)\ar[r]& H^{2d}_c\bigl(X',\Lambda(d)\bigr)\ar[u]_-{\pi'_*}\lefteqn{.}
 }
 \]
 By the commutative diagram 
 \[
  \xymatrix{%
 D'\ar[r]^-{j_0'}\ar[d]_-{\cong}^-{\pi}& X'\ar[d]^-{\pi'}\\
 D\ar[r]^-{j_0}& X\lefteqn{,}
 }
 \]
 we can construct the map
 $\pi_*\colon H^{2d}(D',j_0'^!\Lambda(d))\longrightarrow H^{2d}(D,j_0^!\Lambda(d))$ as the composite
 \begin{align*}
  H^{2d}\bigl(D',j_0'^!\Lambda(d)\bigr)&=H^{2d}\bigl(D,\pi_!j_0'^!\Lambda(d)\bigr)
  =H^{2d}\bigl(D,\pi_!j_0'^!\pi'^!\Lambda(d)\bigr)\\
  &=H^{2d}\bigl(D,\pi_!\pi^!j_0^!\Lambda(d)\bigr)\yrightarrow{\adj}
  H^{2d}\bigl(D,j_0^!\Lambda(d)\bigr).
 \end{align*}
 The commutativity of the following diagram is immediate:
 \[
  \xymatrix{%
 H^{2d}\bigl(D,j_0^!\Lambda(d)\bigr)\ar@{=}[r]&
 H^{2d}_c\bigl(X,Rj_{0!}j_0^!\Lambda(d)\bigr)\ar[r]& H^{2d}_c\bigl(X,\Lambda(d)\bigr)\\
 H^{2d}\bigl(D',j_0'^!\Lambda(d)\bigr)\ar@{=}[r]\ar[u]_-{\pi_*}&
 H^{2d}_c\bigl(X',Rj'_{0!}j_0'^!\Lambda(d)\bigr)\ar[r]& H^{2d}_c\bigl(X',\Lambda(d)\bigr)\ar[u]_-{\pi'_*}\lefteqn{.}
 }
 \]
 Thus it suffices to show that 
 $\pi_*\circ \pi^*\colon H^{2d}(D,j_0^!\Lambda(d))\longrightarrow H^{2d}(D,j_0^!\Lambda(d))$
 is the identity map. This map is induced from the composite of
 $j_0^!\Lambda\yrightarrow{\adj}R\pi_*\pi^*j_0^!\Lambda=\pi_!\pi^!j_0^!\Lambda\yrightarrow{\adj}j_0^!\Lambda$,
 which is the identity map since $R\pi_*=\pi_!$ is the quasi-inverse of $\pi^*=\pi^!$ by the assumption that
 $\pi\colon D'\longrightarrow D$ is an isomorphism.
 Now the proof is complete.
\end{prf}

\begin{prop}\label{prop:Fix-GAGA}
 Let $X$ and $\Gamma$ be purely $d$-dimensional schemes which are separated of finite type over $k$,
 Assume that $X$ (resp.\ $\Gamma$) is smooth (resp.\ generically smooth) over $k$.
 Let $\gamma\colon \Gamma\longrightarrow X\times_kX$ be a morphism over $k$.

 For an open and closed subset $D$ of $\Fix \gamma:=\Gamma\times_{X\times_kX}X$ which is proper over $k$,
 we can define $\#\Fix_D(\gamma)$ in the same way as in Definition \ref{defn:Fix}.
 Then we have $\#\Fix_{D^\ad}(\gamma^\ad)=\#\Fix_D(\gamma)$.
 In particular, if $\Fix\gamma$ is proper over $k$,
 then $\#\Fix(\gamma^\ad)$ coincides with the number of points fixed by $\gamma$ in the usual
 (intersection-theoretic) sense.
\end{prop}

\begin{prf}
 We denote the natural morphism $\Fix\gamma\longrightarrow X$ by $\gamma_0$,
 the open and closed immersion $D\hooklongrightarrow \Fix \gamma$ by $j$, and put $j_0=\gamma_0\circ j$.
 Then the proposition is clear from Proposition \ref{prop:Gysin-GAGA} and the following commutative diagrams:
 \[
  \xymatrix{%
 H^{2d}\bigl(\Gamma,\gamma^!\Lambda(d)\bigr)\ar[r]^-{\delta^*}\ar[d]^-{\varepsilon^*}&
 H^{2d}\bigl(\Fix\gamma,\gamma^!\Lambda(d)\bigr)\ar[r]^-{j^*}\ar[d]^-{\varepsilon^*}&
 H^{2d}\bigl(D,j^*\gamma^!\Lambda(d)\bigr)\ar[d]^-{\varepsilon^*}\\
 H^{2d}\bigl(\Gamma^\ad,\gamma^{\ad!}\Lambda(d)\bigr)\ar[r]^-{\delta^{\ad *}}&
 H^{2d}\bigl(\Fix\gamma^\ad,\gamma^{\ad!}\Lambda(d)\bigr)\ar[r]^-{j^{\ad*}}&
 H^{2d}\bigl(D^\ad,j^{\ad*}\gamma^{\ad!}\Lambda(d)\bigr)\lefteqn{,}
 }
 \]
 \[
  \xymatrix{%
 H^{2d}\bigl(D,j^*\gamma^!\Lambda(d)\bigr)\ar@{=}[r]\ar[d]^-{\varepsilon}&
 H^{2d}_c\bigl(X,Rj_{0!}j_0^!\Lambda(d)\bigr)\ar[r]^-{\adj}\ar[d]^-{\varepsilon^*}&
 H^{2d}_c\bigl(X,\Lambda(d)\bigr)\ar[r]^-{\Tr_X}\ar[d]^-{\varepsilon^*}&\Lambda\ar@{=}[d]\\
 H^{2d}\bigl(D^\ad,j^{\ad*}\gamma^{\ad!}\Lambda(d)\bigr)\ar@{=}[r]&
 H^{2d}_c\bigl(X^\ad,Rj_{0!}^\ad j_0^{\ad!}\Lambda(d)\bigr)\ar[r]^-{\adj}&
 H^{2d}_c\bigl(X^\ad,\Lambda(d)\bigr)\ar[r]^-{\Tr_{X^\ad}}&\Lambda\lefteqn{.}
 }
 \]
\end{prf}

\begin{rem}
 By Proposition \ref{prop:Fix-etale} and Proposition \ref{prop:Fix-GAGA}, we may often calculate $\#\Fix\gamma$.
 For example, we can apply the method in \cite[\S 2.6]{MR2383890} to calculate the number of fixed points
 on some Rapoport-Zink spaces. Since the period space for a Rapoport-Zink space $\mathscr{M}$ is an open adic
 subspace of an algebraic variety, we can use Proposition \ref{prop:Fix-GAGA} for counting fixed points on
 the period space. As the period map from $\mathscr{M}$ to the period space is \'etale,
 Proposition \ref{prop:Fix-etale} enables us to count fixed points on $\mathscr{M}$.
\end{rem}

\begin{defn}
 In the setting introduced at the beginning of this subsection, assume moreover that $\gamma_1$ is proper. %
 We define the action $\gamma^*$ of $\gamma$ on $R\Gamma_c(X,\Lambda)$ as follows:
\[
 \gamma^*\colon R\Gamma_c(X,\Lambda)\yrightarrow{\gamma_1^*}R\Gamma_c(\Gamma,\Lambda)
 \yrightarrow{\Gys_{\gamma_2}} R\Gamma_c(\Gamma,\gamma_2^!\Lambda)=R\Gamma_c(X,R\gamma_{2!}\gamma_2^!\Lambda)
 \yrightarrow{\adj}R\Gamma_c(X,\Lambda).
\]
\end{defn}

\begin{exa}
 Let $f\colon X\longrightarrow X$ be a proper morphism over $S$.
 Then $\gamma_f^*$ (\cf Example \ref{exa:self-mor}) obviously coincides with $f^*$.
\end{exa}

In the sequel we assume that $X$ and $\Gamma$ are quasi-compact and $\gamma_1$ (and hence $\gamma$) is proper.
We will describe $\gamma^*$ by means of a compactification of $\gamma\colon \Gamma\longrightarrow X\times_SX$.

\begin{defn}
 A compactification of $\gamma\colon \Gamma\longrightarrow X\times_SX$ is a triple 
 $(X\hookrightarrow \overline{X},\Gamma\hookrightarrow \overline{\Gamma},\overline{\gamma})$, where
 $X\hooklongrightarrow \overline{X}$ and $\Gamma\hooklongrightarrow \overline{\Gamma}$ are dense open immersions
 into pseudo-adic spaces which are proper over $S$ and $\overline\gamma$ is a proper $S$-morphism
 which makes the following diagram commutative:
 \[
  \xymatrix{%
 \Gamma\ar[r]^-{\gamma}\ar[d]& X\times_SX\ar[d]\\ \overline{\Gamma}\ar[r]^-{\overline{\gamma}}& \overline{X}\times_S \overline{X}\lefteqn{.}
 }
 \]
 For simplicity, we often write $\overline{\gamma}\colon \overline{\Gamma}\longrightarrow \overline{X}\times_S\overline{X}$ for
 $(X\hookrightarrow \overline{X},\Gamma\hookrightarrow \overline{\Gamma},\overline{\gamma})$.
\end{defn}

\begin{rem}
 For a compactification 
 $\overline{\gamma}\colon \overline{\Gamma}\longrightarrow \overline{X}\times_S\overline{X}$,
 we have $\Gamma=\overline{\gamma}_1^{\,-1}(X)$. Indeed, since $\gamma_1$ is proper, 
 the open immersion $\Gamma\hooklongrightarrow \overline{\gamma}_1^{\,-1}(X)$ is proper.
 On the other hand, since $\Gamma$ is assumed to be dense in $\overline{\Gamma}$, it is also dense
 in $\overline{\gamma}_1^{\,-1}(X)$. Thus we have $\Gamma=\overline{\gamma}_1^{\,-1}(X)$.
 In other words, $\overline{\gamma}_1^{\,-1}(X)$ is contained in $\overline{\gamma}_2^{\,-1}(X)$.
\end{rem}

\begin{exa}
 Let $X\hooklongrightarrow X^c$ and $\Gamma\hooklongrightarrow \Gamma^c$ be the universal compactifications
 of $X$ and $\Gamma$ over $S$, respectively
 (\cf \cite[Definition 5.1.1, Theorem 5.1.5, Corollary 5.1.6]{MR1734903}).
 Then the morphism $\gamma^c\colon \Gamma^c\longrightarrow X^c\times_SX^c$ over $S$ is
 naturally induced from $\gamma$,
 and gives a compactification of $\gamma\colon \Gamma\longrightarrow X\times_SX$.

 If $\gamma\colon \Gamma\longrightarrow X\times_SX$ can be extended to a morphism 
 $\gamma'\colon \Gamma'\longrightarrow X'\times_SX'$ where $X'$ (resp.\ $\Gamma'$) is an adic space which is
 partially proper taut over $S$ and contains $X$ (resp.\ $\Gamma$) as an open adic subspace, 
 then we can construct another compactification of $\gamma$. Let $\overline{X}$ (resp.\ $\overline{\Gamma}$)
 be the closure of $X$ (resp.\ $\Gamma$) in $X'$ (resp.\ $\Gamma'$) and regard it as a pseudo-adic space.
 Then $\overline{X}$ and $\overline{\Gamma}$ are proper over $S$ and $\gamma'$ induces a morphism
 $\overline{\gamma}\colon \overline{\Gamma}\longrightarrow \overline{X}\times_S\overline{X}$.
 It gives a compactification of $\gamma$. This construction should be more convenient for practical use.
\end{exa}

Take a compactification 
$\overline{\gamma}\colon \overline{\Gamma}\longrightarrow \overline{X}\times_S\overline{X}$ of $\gamma$ and
denote the open immersion $X\hooklongrightarrow \overline{X}$ by $j$.
We denote the open immersions $X\times_SX\hooklongrightarrow \overline{X}\times_SX$,
$X\times_S\overline{X}\hooklongrightarrow \overline{X}\times_S\overline{X}$ by $j\times 1$,
and $X\times_SX\hooklongrightarrow X\times_S\overline{X}$,
$\overline{X}\times_SX\hooklongrightarrow \overline{X}\times_S\overline{X}$ by $1\times j$.
Consider the natural isomorphisms
\begin{align*}
 H^{2d}_{\overline{\gamma}(\overline{\Gamma})}\bigl(\overline{X}\times_S\overline{X},(j\times 1)_*(1\times j)_!\Lambda(d)\bigr)
 &\yrightarrow{\cong}H^{2d}_{\gamma(\Gamma)}\bigl(X\times_S\overline{X},(1\times j)_!\Lambda(d)\bigr)\\
 &\yrightarrow{\cong}H^{2d}_{\gamma(\Gamma)}\bigl(X\times_SX,\Lambda(d)\bigr).
\end{align*}
Note that the first isomorphy is a consequence of 
$\Gamma=\overline{\gamma}_1^{\,-1}(X)=\overline{\gamma}^{\,-1}(X\times_S\overline{X})$.
We also denote by $\cl(\gamma)$ the element of $H^{2d}_{\overline{\gamma}(\overline{\Gamma})}(\overline{X}\times_S\overline{X},(j\times 1)_*(1\times j)_!\Lambda(d))$ that is mapped to $\cl(\gamma)$
by the homomorphism above.
Since the projection formula gives
\begin{align*}
 (j\times 1)_!\Lambda\Lotimes (j\times 1)_*(1\times j)_!\Lambda
 &=(j\times 1)_!\bigl(\Lambda\Lotimes (j\times 1)^*(j\times 1)_*(1\times j)_!\Lambda\bigr)\\
 &=(j\times j)_!\Lambda,
\end{align*}
the cup product with $\cl(\gamma)$ induces the map
\[
 R\Gamma\bigl(\overline{X}\times_S\overline{X},(j\times 1)_!\Lambda\bigr)\longrightarrow R\Gamma\bigl(\overline{X}\times_S\overline{X},(j\times j)_!\Lambda(d)[2d]\bigr).
\]

\begin{prop}\label{prop:corr-cl}
 The map $\gamma^*$ coincides with the composite below:
 \begin{align*}
  R\Gamma_c(X,\Lambda)&=R\Gamma(\overline{X},j_!\Lambda)\yrightarrow{\pr_1^*}R\Gamma(\overline{X}\times_S\overline{X},\pr_1^*j_!\Lambda)
 =R\Gamma\bigl(\overline{X}\times_S\overline{X},(j\times 1)_!\Lambda\bigr)\\
  &\yrightarrow{\cup \cl(\gamma)}R\Gamma\bigl(\overline{X}\times_S\overline{X},(j\times j)_!\Lambda(d)[2d]\bigr)
  =R\Gamma_c\bigl(X\times_SX,\Lambda(d)[2d]\bigr)\\
  &\yrightarrow{\Gys_{\pr_2}}R\Gamma_c(X\times_SX,\pr_2^!\Lambda)=R\Gamma_c(X,R\pr_{2!}\pr_2^!\Lambda)
  \yrightarrow{\adj}R\Gamma_c(X,\Lambda).
 \end{align*}
\end{prop}

\begin{prf}
 We also denote the open immersion $\Gamma\hooklongrightarrow \overline{\Gamma}$ by $j$.
 First let us prove the commutativity of the following diagram:
\[
 \xymatrix{%
 R\Gamma\bigl(\overline{X}\times_S\overline{X},(j\times 1)_!\Lambda\bigr)\ar[r]^-{\cup\cl(\gamma)}\ar[d]^-{\overline{\gamma}^*}& R\Gamma\bigl(\overline{X}\times_S\overline{X},(j\times j)_!\Lambda(d)[2d]\bigr)\\
 R\Gamma\bigl(\overline{\Gamma},\overline{\gamma}^*(j\times 1)_!\Lambda\bigr)\ar@{=}[d]&
 R\Gamma\bigl(\overline{X}\times_S\overline{X},(j\times j)_!R\gamma_!\gamma^!\Lambda(d)[2d]\bigr)
 \ar@{=}[d]\ar[u]_-{\adj}\\
 R\Gamma(\overline{\Gamma},j_!\Lambda)\ar[r]^-{j_!\Gys_{\gamma}}& R\Gamma\bigl(\overline{\Gamma},j_!\gamma^!\Lambda(d)[2d]\bigr)\lefteqn{.}
 }
\]
 By the adjointness of $(j\times 1)_!$ and $(j\times 1)^*$,
 it suffices to show the commutativity of the following diagram, where $\gamma':=(1\times j)\circ \gamma$:
 \[
  \xymatrix{%
 \Lambda\ar[r]^-{\cl(\gamma)}\ar[d]^-{\adj}&(1\times j)_!\Lambda(d)[2d]\\
 R\gamma'_*\gamma'^*\Lambda\ar@{=}[d]& (1\times j)_!R\gamma_!\gamma^!\Lambda(d)[2d]\ar@{=}[d]\ar[u]_-{\adj}\\
 R\gamma'_*\Lambda\ar[r]^-{R\gamma'_*\Gys_{\gamma}}& R\gamma'_*\gamma^!\Lambda(d)[2d]\lefteqn{.}
 }
 \]
 Here $\cl(\gamma)$ is regarded as an element of $\Hom(\Lambda,(1\times j)_!\Lambda(d)[2d])$
 by the maps
 \begin{align*}
  &H^{2d}_{\gamma(\Gamma)}(X\times_SX,\Lambda(d))\yleftarrow{\cong}H^{2d}_{\gamma(\Gamma)}(X\times_S\overline{X},(1\times j)_!\Lambda(d))\\
  &\qquad\longrightarrow H^{2d}(X\times_S\overline{X},(1\times j)_!\Lambda(d))=\Hom\bigl(\Lambda,(1\times j)_!\Lambda(d)[2d]\bigr).
 \end{align*} 
 By the construction, it is obtained by the composite
 \begin{align*}
  \Lambda&\yrightarrow{\adj} R\gamma'_*\gamma'^*\Lambda=R\gamma'_*\Lambda\yrightarrow{R\gamma'_*\Gys_{\gamma}}
  R\gamma'_*\gamma^!\Lambda(d)[2d]=R\gamma'_!\gamma'^!(1\times j)_!\Lambda(d)[2d]\\
  &\yrightarrow{\adj}(1\times j)_!\Lambda(d)[2d].
 \end{align*}
 Since it is easy to see that two maps
 \begin{gather*}
  R\gamma'_*\gamma^!\Lambda=R\gamma'_!\gamma'^!(1\times j)_!\Lambda\yrightarrow{\adj}(1\times j)_!\Lambda,\\
  R\gamma'_*\gamma^!\Lambda=(1\times j)_!R\gamma_!\gamma^!\Lambda\yrightarrow{\adj}(1\times j)_!\Lambda
 \end{gather*}
 coincide, we have the desired commutativity.

 Therefore, the composite of
 \begin{align*}
  R\Gamma_c(X,\Lambda)&=R\Gamma(\overline{X},j_!\Lambda)\yrightarrow{\pr_1^*}R\Gamma(\overline{X}\times_S\overline{X},\pr_1^*j_!\Lambda)
 =R\Gamma\bigl(\overline{X}\times_S\overline{X},(j\times 1)_!\Lambda\bigr)\\
  &\yrightarrow{\cup \cl(\gamma)}R\Gamma\bigl(\overline{X}\times_S\overline{X},(j\times j)_!\Lambda(d)[2d]\bigr)
  =R\Gamma_c\bigl(X\times_SX,\Lambda(d)[2d]\bigr)
 \end{align*}
 coincides with the composite of
  \begin{align*}
  R\Gamma_c(X,\Lambda)&\yrightarrow{\gamma_1^*} R\Gamma_c(\Gamma,\Lambda)
   \yrightarrow{\Gys_\gamma}R\Gamma_c\bigl(\Gamma,\gamma^!\Lambda(d)[2d]\bigr)
   =R\Gamma_c\bigl(X\times_SX,R\gamma_!\gamma^!\Lambda(d)[2d]\bigr)\\
   &\yrightarrow{\adj}R\Gamma_c\bigl(X\times_SX,\Lambda(d)[2d]\bigr).
 \end{align*}
 Therefore it suffices to show that the composite of
 \begin{align*}
  R\Gamma_c(\Gamma,\Lambda)&\yrightarrow{\Gys_\gamma}R\Gamma_c\bigl(\Gamma,\gamma^!\Lambda(d)[2d]\bigr)=R\Gamma_c\bigl(X\times_SX,R\gamma_!\gamma^!\Lambda(d)[2d]\bigr)\\
  &\yrightarrow{\adj}R\Gamma_c\bigl(X\times_SX,\Lambda(d)[2d]\bigr)\yrightarrow{\Gys_{\pr_2}}
  R\Gamma_c\bigl(X\times_SX,\pr_2^!\Lambda\bigr)\\
  &=R\Gamma_c\bigl(X,R\pr_{2!}\pr_2^!\Lambda\bigr)\yrightarrow{\adj}R\Gamma_c(X,\Lambda)
 \end{align*}
 is equal to the composite of
 \[
  R\Gamma_c(\Gamma,\Lambda)\yrightarrow{\Gys_{\gamma_2}}R\Gamma_c(\Gamma,\gamma_2^!\Lambda)
 =R\Gamma_c(X,R\gamma_{2!}\gamma_2^!\Lambda)\yrightarrow{\adj}R\Gamma_c(X,\Lambda).
 \]
 It is an easy consequence of $\Gys_{\gamma_2}=\gamma^!\Gys_{\pr_2}(d)[2d]\circ \Gys_{\gamma}$, which
 can be proved directly from the construction of the Gysin maps.
\end{prf}

\section{Lefschetz trace formula for open adic spaces}
\subsection{K\"unneth formula}
\begin{lem}\label{lem:constant}
 Let $X$ be a finite-dimensional pseudo-adic space which is quasi-separated of weakly finite type over $S$. 
 Let $\mathcal{F}$ be a $\Lambda$-sheaf on $X$ and $L$ a bounded complex of $\Lambda$-modules.
 Then we have an isomorphism $R\Gamma(X,\mathcal{F}\Lotimes L_X)\cong R\Gamma(X,\mathcal{F})\Lotimes L$.
\end{lem}

\begin{prf}
 We may assume $L$ is a $\Lambda$-module.
 By \cite[Corollary 2.8.3, Corollary 1.8.8]{MR1734903},
 the cohomological dimension of $R\Gamma$ is finite. Therefore, by taking a free resolution of $L$,
 we may reduce to the case where $L$ is a free $\Lambda$-module.
 Then the claim holds, since $R\Gamma$ commutes with any direct sum if $X$ is quasi-compact and quasi-separated
 (\cite[Lemma 2.3.13 i)]{MR1734903}). 
\end{prf}

Similarly, we can prove the following:

\begin{lem}\label{lem:constant-cpt-supp}
 Let $X$ be a finite-dimensional pseudo-adic space which is separated, locally of +-weakly finite type and taut over $S$. 
 Let $\mathcal{F}$ be a $\Lambda$-sheaf on $X$ and $L$ a bounded complex of $\Lambda$-modules. 
 Then we have an isomorphism $R\Gamma_c(X,\mathcal{F}\Lotimes L_X)\cong R\Gamma_c(X,\mathcal{F})\Lotimes L$.
\end{lem}

\begin{prf}
 Since $R\Gamma_c$ has finite cohomological dimension (\cite[Proposition 5.5.8, Corollary 1.8.8]{MR1734903})
 and commutes with any direct sum (\cite[Proposition 5.4.5 i)]{MR1734903}),
 the proof is exactly the same as the previous lemma.
\end{prf}

\begin{cor}\label{cor:perfect-cpx}
 Let $X$ be a finite-dimensional pseudo-adic space which is separated, locally of +-weakly finite type and taut over $S$.
 Assume that $H^i_c(X,\Lambda)$ is a finitely generated $\Lambda$-module for every $i$.
 Then $R\Gamma_c(X,\Lambda)$ is a perfect $\Lambda$-complex.

 In particular, if $\underline{X}$ is separated of finite type over $S$ and $\lvert X\rvert$ is
 a locally closed constructible subset of $\underline{X}$, then $R\Gamma_c(X,\Lambda)$ is a perfect $\Lambda$-complex.
\end{cor}

\begin{prf}
 Since $R\Gamma_c(X,\Lambda)$ is bounded with finitely generated cohomology,
 by \cite[{[Rapport]}, Lemme 4.5.1]{MR0463174} it suffices to show that $R\Gamma_c(X,\Lambda)$ has finite tor-dimension. 
 Let $M$ be a $\Lambda$-module. Then Lemma \ref{lem:constant-cpt-supp} says
 $R\Gamma_c(X,\Lambda)\Lotimes M\cong R\Gamma_c(X,M_X)$.
 In particular, the $i$th cohomology of $R\Gamma_c(X,\Lambda)\Lotimes M$ vanishes unless $0\le i\le 2\dim X$.
 This means that $R\Gamma_c(X,\Lambda)$ has finite tor-dimension.

 The latter part follows from the finiteness results due to Huber
 (\cite[Corollary 2.3]{MR1620114}, \cite[Corollary 5.4]{MR2336836}).
\end{prf}

Let $X$ and $Y$ be finite-dimensional pseudo-adic spaces which are quasi-separated of weakly finite type over $S$,
and denote their structure maps by $a\colon X\longrightarrow S$ and $b\colon Y\longrightarrow S$.
Denote the first (resp.\ second) projection by $\pr_1\colon X\times_SY\longrightarrow X$ (resp.\ $\pr_2\colon X\times_SY\longrightarrow Y$).
Let $\mathcal{F}$ be a sheaf on $X$ and $\mathcal{G}$ a sheaf on $Y$. 
Put $\mathcal{F}\Lboxtimes\mathcal{G}=\pr_1^*\mathcal{F}\Lotimes\pr_2^*\mathcal{G}$.
Then we have the canonical homomorphism 
\[
 R\Gamma(X,\mathcal{F})\Lotimes R\Gamma(Y,\mathcal{G})\longrightarrow R\Gamma(X\times_SY,\mathcal{F}\Lboxtimes\mathcal{G}),
\]
which is called the K\"unneth homomorphism.

\begin{lem}\label{lem:criterion-Kunneth-isom}
 If the canonical map
 \[
  \mathcal{F}\Lotimes R\pr_{1*}\pr_2^*\mathcal{G} \longrightarrow R\pr_{1*}(\pr_1^*\mathcal{F}\Lotimes \pr_2^*\mathcal{G})
 \]
 is an isomorphism, the K\"unneth homomorphism is also an isomorphism.
\end{lem}

\begin{prf}
 By the quasi-compact/generalizing base change theorem (\cite[Theorem 4.3.1]{MR1734903}), we have 
 $R\pr_{1*}\pr_2^*\mathcal{G}\cong R\Gamma(Y,\mathcal{G})_X$.
 Then by Lemma \ref{lem:constant}, we have an isomorphism
 \begin{align*}
  R\Gamma(X,\mathcal{F})\Lotimes R\Gamma(Y,\mathcal{G})&\cong R\Gamma\bigl(X,\mathcal{F}\Lotimes R\Gamma(Y,\mathcal{G})_X\bigr)\cong R\Gamma\bigl(X,\mathcal{F}\Lotimes R\pr_{1*}\pr_2^*\mathcal{G}\bigr)\\
  &\yrightarrow{\cong} R\Gamma\bigl(X,R\pr_{1*}(\pr_1^*\mathcal{F}\Lotimes \pr_2^*\mathcal{G})\bigr)
  \cong R\Gamma(X\times_SY,\mathcal{F}\Lboxtimes \mathcal{G}).
 \end{align*}
 It is easy to see that the isomorphism above is actually the K\"unneth homomorphism.
\end{prf}

\begin{prop}\label{prop:Kunneth-const}
 In the case $\mathcal{F}=\Lambda$, the K\"unneth homomorphism
 \[
 R\Gamma(X,\Lambda)\Lotimes R\Gamma(Y,\mathcal{G})\longrightarrow R\Gamma(X\times_SY,\Lambda\Lboxtimes\mathcal{G})
 \]
 is an isomorphism.
\end{prop}

\begin{prf}
 Clear by Lemma \ref{lem:criterion-Kunneth-isom}.
\end{prf}

\begin{prop}\label{prop:Kunneth-2}
 Let $j\colon U\hooklongrightarrow Y$ be a quasi-compact open immersion.
 If $\mathcal{G}=j_!\Lambda$, then the K\"unneth homomorphism
 \[
 R\Gamma(X,\mathcal{F})\Lotimes R\Gamma(Y,j_!\Lambda)\longrightarrow R\Gamma(X\times_SY,\mathcal{F}\Lboxtimes j_!\Lambda)
 \]
 is an isomorphism.
\end{prop}

\begin{prf}
 By Lemma \ref{lem:criterion-Kunneth-isom}, we have only to prove that the natural map
 \[
 (R\pr_{2*}\pr_1^*\mathcal{F})\Lotimes j_!\Lambda\longrightarrow R\pr_{2*}(\pr_1^*\mathcal{F}\Lotimes \pr_2^*j_!\Lambda)
 \]
 is an isomorphism. 
 Put $Z=(\underline{Y},\lvert Y\rvert\setminus \lvert U\rvert)$ and denote the closed immersion
 $Z\hooklongrightarrow Y$ of pseudo-adic spaces by $i$. 
 By the distinguished triangle $j_!\Lambda\longrightarrow \Lambda\longrightarrow i_*\Lambda\yrightarrow{+1} j_!\Lambda[1]$,
 it suffices to prove that the natural map
 \[
 (R\pr_{2*}\pr_1^*\mathcal{F})\Lotimes i_*\Lambda\longrightarrow R\pr_{2*}(\pr_1^*\mathcal{F}\Lotimes \pr_2^*i_*\Lambda)
 \]
 is an isomorphism.

 Consider the following commutative diagram whose rectangles are cartesian:
 \[
  \xymatrix{%
 X\times_SZ\ar[r]^-{1\times i}\ar[d]^-{\pr'_2}& X\times_SY\ar[r]^-{\pr_1}\ar[d]^-{\pr_2}& X\ar[d]^-{a}\\
 Z\ar[r]^-{i}& Y\ar[r]^-{b}& S\lefteqn{.}
 }
 \]
 By the quasi-compact/generalizing base change theorem, the base change map
 $i^*R\pr_{2*}\pr_1^*\mathcal{F}\longrightarrow i^*b^*Ra_*\mathcal{F}$
 and $R\pr'_{2*}(1\times i)^*\pr_1^*\mathcal{F}\longrightarrow i^*b^*Ra_*\mathcal{F}$ are
 isomorphisms. Thus the base change map
 $i^*R\pr_{2*}\pr_1^*\mathcal{F}\longrightarrow R\pr'_{2*}(1\times i)^*\pr_1^*\mathcal{F}$
 is also an isomorphism. By this, we have
 \begin{align*}
  (R\pr_{2*}\pr_1^*\mathcal{F})\Lotimes i_*\Lambda&\cong i_*i^*R\pr_{2*}\pr_1^*\mathcal{F}\yrightarrow{\cong}
  i_*R\pr'_{2*}(1\times i)^*\pr_1^*\mathcal{F}\\
  &=R\pr_{2*}(1\times i)_*(1\times i)^*\pr_1^*\mathcal{F}
  \cong R\pr_{2*}(\pr_1^*\mathcal{F}\Lotimes (1\times i)_*\Lambda)\\
  &\cong R\pr_{2*}(\pr_1^*\mathcal{F}\Lotimes \pr_2^*i_*\Lambda),
 \end{align*}
 which completes the proof.
\end{prf}

\begin{cor}\label{cor:Kunneth-square-shape}
 Let $U\subset X$ and $V\subset Y$ be quasi-compact open adic subspaces.
 Denote the open immersions $U\hooklongrightarrow X$, $V\hooklongrightarrow Y$ by $j$, $j'$ respectively.
 We write $j\times 1$ for $U\times_S V\hooklongrightarrow X\times_S V$ and $U\times_S Y\hooklongrightarrow X\times_S Y$,
 $1\times j'$ for $U\times_S V\hooklongrightarrow U\times_S Y$ and $X\times_S V\hooklongrightarrow X\times_S Y$.
 Since $(j\times 1)^*(j_*\Lambda\Lboxtimes j'_!\Lambda)\cong (j\times 1)^*\pr_2^*j'_!\Lambda\cong (1\times j')_!\Lambda$, we have the canonical morphism 
 $j_*\Lambda\Lboxtimes j'_!\Lambda\longrightarrow (j\times 1)_*(1\times j')_!\Lambda$ which is denoted by $\tau$.
 Then the map
 \[
 R\Gamma(X\times_SY,j_*\Lambda\Lboxtimes j'_!\Lambda)\longrightarrow R\Gamma\bigl(X\times_SY,(j\times 1)_*(1\times j')_!\Lambda\bigr)
 \]
 induced by $\tau$ is an isomorphism.
\end{cor}

\begin{prf}
 Clear from Proposition \ref{prop:Kunneth-2} and the commutative diagram below: 
 \[
 \xymatrix{%
 R\Gamma\bigl(X\times_SY,j_*\Lambda\Lboxtimes j'_!\Lambda\bigr)\ar[r]&R\Gamma\bigl(X\times_SY,(j\times 1)_*(1\times j')_!\Lambda\bigr)\ar[d]^-{\cong}\\
 & R\Gamma\bigl(U\times_SY,(1\times j')_!\Lambda\bigr)\\
 R\Gamma(X,j_*\Lambda)\Lotimes R\Gamma(Y,j'_!\Lambda)\ar[uu]_{\cong}\ar[r]^-{\cong}& R\Gamma(U,\Lambda)\Lotimes R\Gamma(Y,j'_!\Lambda)\ar[u]_-{\cong}\lefteqn{.}
 }
 \]
\end{prf}

\begin{rem}\label{rem:isom-failure}
 Unlike the case of schemes, the morphism $\tau$ itself is not an isomorphism in general.
 For example, put $U=V=\mathbb{D}^1$ and let $X$, $Y$ be the closure $\overline{\mathbb{D}^1}$ of $\mathbb{D}^1$ in 
 $(\mathbb{A}^1)^{\mathrm{ad}}$. Take a continuous valuation $\lvert\ \rvert\colon k\longrightarrow \R_{\ge 0}$
 on $k$ and consider the point $x$ of $(\mathbb{A}^1)^{\mathrm{ad}}$ corresponding to the valuation
 \[
  k[T]\longrightarrow \R_{\ge 0}\times \Z;\quad \sum_i a_iT^i\longmapsto \max\bigl\{(\lvert a_i\rvert,i)\bigr\}
 \]
 (here we endow $\R_{\ge 0}\times \Z$ with the lexicographic order).
 Then $x$ lies in the closure of $\mathbb{D}^1$ in $(\A^1)^{\mathrm{ad}}$ but does not lie in $\mathbb{D}^1$ itself.
 In $(\A^1)^{\mathrm{ad}}$ it has a unique generalization $y$ which is given by the following valuation: 
 \[
  k[T]\longrightarrow \R_{\ge 0};\quad \sum_i a_iT^i\longmapsto \max\{\lvert a_i\rvert\}.
 \]
 It is easy to see that $y$ belongs to $\mathbb{D}^1$.

 By the diagonal map $(\A^1)^{\mathrm{ad}}\longrightarrow (\A^2)^{\mathrm{ad}}$, we regard $x$ and $y$ as points of
 $(\A^2)^{\mathrm{ad}}$. Then $x$ lies in $\overline{\mathbb{D}^1}\times_S\overline{\mathbb{D}^1}
 =\pr_1^{-1}(\overline{\mathbb{D}^1})\cap \pr_2^{-1}(\overline{\mathbb{D}^1})$ and
 $y$ is a unique generalization of $x$ in $(\A^2)^{\mathrm{ad}}$. 
 Let $j$ be the open immersion $\mathbb{D}^1\hooklongrightarrow \overline{\mathbb{D}^1}$, which is quasi-compact.
 It is clear that $(j_*\Lambda\Lboxtimes j'_!\Lambda)_{\overline{x}}=0$. 
 On the other hand, by \cite[Proposition 2.6.4]{MR1734903} and its proof,
 we have $((j\times 1)_*(1\times j)_!\Lambda)_{\overline{x}}=((1\times j)_!\Lambda)_{\overline{y}}=\Lambda$.
 Thus $\tau$ is not an isomorphism.
\end{rem}
 
\subsection{Unlocalized Lefschetz trace formula}
 In the remaining part of this section, we use the same notation as in \S \ref{subsec:corr};
 let $\gamma\colon \Gamma\longrightarrow X\times_SX$ be an $S$-morphism between purely $d$-dimensional
 adic spaces which are separated of finite type over $S$, and assume that $X$ (resp.\ $\Gamma$) 
 is smooth (resp.\ generically smooth) over $S$ and $\gamma_1$ is proper.
 Fix a compactification $\overline{\gamma}\colon \overline{\Gamma}\longrightarrow \overline{X}\times_S\overline{X}$ and denote the open immersion $X\hooklongrightarrow \overline{X}$ by $j$.

 As in \S 2, $\gamma$ defines the element $\cl(\gamma)$ of $H^{2d}(\overline{X}\times_S\overline{X},(j\times 1)_*(1\times j)_!\Lambda(d))$. 
 Moreover, by Corollary \ref{cor:Kunneth-square-shape},
 the map $H^{2d}(\overline{X}\times_S\overline{X},j_*\Lambda \Lboxtimes j_!\Lambda(d))\longrightarrow H^{2d}(\overline{X}\times_S\overline{X},(j\times 1)_*(1\times j)_!\Lambda(d))$ induced by the canonical morphism $\tau\colon j_*\Lambda \Lboxtimes j_!\Lambda\longrightarrow (j\times 1)_*(1\times j)_!\Lambda$ is an isomorphism. 
 Therefore $\gamma$ defines an element of
 $H^{2d}(Y\times_SY,(1\times j)_!(j\times 1)_*\Lambda(d))$ denoted by $[\gamma]$.
 The diagonal morphism $\overline{\delta}\colon \overline{X}\longrightarrow \overline{X}\times_S\overline{X}$
 induces the pull-back map 
 $\overline{\delta}^*\colon H^{2d}(\overline{X}\times_S\overline{X},j_*\Lambda \Lboxtimes j_!\Lambda(d))\longrightarrow H^{2d}(\overline{X},j_!\Lambda(d))=H^{2d}_c(X,\Lambda(d))$.

\begin{prop}\label{prop:LTF-unlocalized}
 In the situation above, we have the equality
 \[
 \Tr\bigl(\gamma^*;R\Gamma_c(X,\Lambda)\bigr)=\Tr_X\bigl(\overline{\delta}^*[\gamma]\bigr).
 \]
\end{prop}

The proof of this proposition is similar to that in the scheme case. However, we will 
include it for the completeness. Let us consider the following diagram:
\[
 \xymatrix{%
 R\Gamma\bigl(X,\Lambda(d)[2d]\bigr)\Lotimes R\Gamma_c(X,\Lambda)\ar[d]^-{\cong}
 \ar@<-2pt>[r]^-{\cong}& R\Hom\bigl(R\Gamma_c(X,\Lambda),\Lambda\bigr)\Lotimes R\Gamma_c(X,\Lambda)\ar[dd]\\
 R\Gamma\bigl(\overline{X}\times_S\overline{X},j_*\Lambda\Lboxtimes j_!\Lambda(d)[2d]\bigr)\ar[d]\ar@/_9pc/[dd]_-{\overline{\delta}^*}\\
 R\Gamma\bigl(\overline{X}\times_S\overline{X},(j\times 1)_*(1\times j)_!\Lambda(d)[2d]\bigr)\ar[r]^-{(*)}& R\Hom\bigl(R\Gamma_c(X,\Lambda),
 R\Gamma_c(X,\Lambda)\bigr)\ar[d]^-{\Tr}\\
 R\Gamma_c\bigl(X,\Lambda(d)[2d]\bigr)\ar[r]^-{\Tr_X}& \Lambda\lefteqn{.}
 }
\]
Here $(*)$ is the map induced by 
\begin{align*}
 &R\Gamma_c(X,\Lambda)\Lotimes R\Gamma\bigl(\overline{X}\times_S\overline{X},(j\times 1)_*(1\times j)_!\Lambda(d)[2d]\bigr)\\
 &\qquad \yrightarrow{\pr_1^*\otimes \id} R\Gamma\bigl(\overline{X}\times_S\overline{X},(j\times 1)_!\Lambda\bigr)\Lotimes R\Gamma\bigl(\overline{X}\times_S\overline{X},(j\times 1)_*(1\times j)_!\Lambda(d)[2d]\bigr)\\
 &\qquad \yrightarrow{\cup}R\Gamma\bigl(\overline{X}\times_S\overline{X},(j\times 1)_!(1\times j)_!\Lambda(d)[2d]\bigr)
 =R\Gamma_c\bigl(X\times_SX,\Lambda(d)[2d]\bigr)\\
 &\qquad \yrightarrow{\pr_{2*}}R\Gamma_c\bigl(X,\Lambda\bigr),
\end{align*}
where $\pr_{2*}$ is the composite of
\begin{align*}
 R\Gamma_c\bigl(X\times_SX,\Lambda(d)[2d]\bigr)&\yrightarrow{\Gys_{\pr_2}}
 R\Gamma_c\bigl(X\times_SX,\pr_2^!\Lambda\bigr)=R\Gamma_c(X,R\pr_{2!}\pr_2^!\Lambda)\\
 &\yrightarrow{\adj}R\Gamma_c(X,\Lambda),
\end{align*}
or equivalently, the composite of
\[
 R\Gamma_c\bigl(X\times_SX,\Lambda(d)[2d]\bigr)=
 R\Gamma_c\bigl(X,R\pr_{2!}\Lambda(d)[2d]\bigr)\yrightarrow{\Tr_{\pr_2}}R\Gamma_c(X,\Lambda).
\]
Since $R\Gamma_c(X,\Lambda)$ is a perfect $\Lambda$-complex
(Corollary \ref{cor:perfect-cpx}), we may define the map 
\[
 \Tr\colon R\Hom(R\Gamma_c(X,\Lambda),R\Gamma_c(X,\Lambda))\longrightarrow \Lambda
\]
(\cite[Expose I]{SGA6}).
It is a unique map such that the composite
\[
 R\Hom\bigl(R\Gamma_c(X,\Lambda),\Lambda\bigr)\Lotimes R\Gamma_c(X,\Lambda)\yrightarrow{\cong} 
 R\Hom\bigl(R\Gamma_c(X,\Lambda),R\Gamma_c(X,\Lambda)\bigr)\yrightarrow{\Tr}\Lambda
\]
coincides with the evaluation map $\mathrm{ev}$.

By Proposition \ref{prop:corr-cl}, $(*)$ maps $\cl(\gamma)$ to $\gamma^*$.
Thus the proposition follows from the commutativity of the lower part of the diagram above,
i.e., the commutativity of the diagram below:
\[
 \xymatrix{%
 R\Gamma\bigl(\overline{X}\times_S\overline{X},j_*\Lambda \Lboxtimes j_!\Lambda(d)[2d]\bigr)\ar[d]\ar[r]^-{\overline{\delta}^*}& R\Gamma_c\bigl(X,\Lambda(d)[2d]\bigr)\ar[dd]^-{\Tr_X}\\
 R\Gamma\bigl(\overline{X}\times_S\overline{X},(j\times 1)_*(1\times j)_!\Lambda(d)[2d]\bigr)\ar[d]^-{(*)}\\
 R\Hom\bigl(R\Gamma_c(X,\Lambda),R\Gamma_c(X,\Lambda)\bigr)\ar[r]^-{\Tr}& \Lambda\lefteqn{.}
 }
\]
To show it, it is sufficient to show the commutativities of the upper part and the outer part of
the diagram above. We will divide their proofs into two propositions:

\begin{prop}\label{prop:commutativity-part-1}
 The following diagram is commutative:
 \[
  \xymatrix{%
  R\Gamma\bigl(X,\Lambda(d)[2d]\bigr)\Lotimes R\Gamma_c(X,\Lambda)\ar[d]^-{\cong}
 \ar@<-2pt>[r]^-{\cong}& R\Hom\bigl(R\Gamma_c(X,\Lambda),\Lambda\bigr)\Lotimes R\Gamma_c(X,\Lambda)\ar[dd]^-{\mathrm{ev}}\\
 R\Gamma\bigl(\overline{X}\times_S\overline{X},j_*\Lambda\Lboxtimes j_!\Lambda(d)[2d]\bigr)\ar[d]^-{\overline{\delta}^*}\\
 R\Gamma_c\bigl(X,\Lambda(d)[2d]\bigr)\ar[r]^-{\Tr_X}& \Lambda\lefteqn{.}
 }
 \]
\end{prop}

\begin{prf}
 Notice that the composite of the left vertical arrows is nothing but the cup product.
 Thus the lemma follows from the next lemma, which is easy to see.
\end{prf}

\begin{lem}
 Let $K$ (resp.\ $L$) be an object of $D^+(\Lambda)$ (resp.\ $D^-(\Lambda)$) and
 $\Phi\colon K\Lotimes L\longrightarrow \Lambda$ a map. We have a natural map 
 $\Phi'\colon K\longrightarrow R\Hom(L,\Lambda)$
 induced from $\Phi$. Then the following diagram is commutative:
 \[
  \xymatrix{%
 K\Lotimes L\ar@<-3pt>[rr]^-{\Phi'\otimes \id}\ar[rd]_-{\Phi}& & R\Hom(L,\Lambda)\Lotimes L\ar[ld]^-{\mathrm{ev}}\\
 &\Lambda \lefteqn{.}
 }
 \]
\end{lem}

\begin{prop}\label{prop:commutativity-part-2}
 The following diagram is commutative:
 \[
 \xymatrix{%
 R\Gamma\bigl(X,\Lambda(d)[2d]\bigr)\Lotimes R\Gamma_c(X,\Lambda)\ar[d]^-{\cong}
 \ar@<-2pt>[r]^-{\cong}& R\Hom\bigl(R\Gamma_c(X,\Lambda),\Lambda\bigr)\Lotimes R\Gamma_c(X,\Lambda)\ar[dd]\\
 R\Gamma\bigl(\overline{X}\times_S\overline{X},j_*\Lambda\Lboxtimes j_!\Lambda(d)[2d]\bigr)\ar[d]\\
 R\Gamma\bigl(\overline{X}\times_S\overline{X},(j\times 1)_*(1\times j)_!\Lambda(d)[2d]\bigr)\ar[r]^-{(*)}& R\Hom\bigl(R\Gamma_c(X,\Lambda),
 R\Gamma_c(X,\Lambda)\bigr)\lefteqn{.}
 }
\]
\end{prop}

\begin{prf}
 Recall that by the definition the morphism $(*)$ is decomposed as
 \begin{align*}
  R\Gamma\bigl(\overline{X}\times_S\overline{X},(j\times 1)_*(1\times j)_!\Lambda(d)[2d]\bigr)&\longrightarrow 
  R\Hom\bigl(R\Gamma_c(X,\Lambda),R\Gamma_c(X\times_SX,\Lambda(d)[2d])\bigr)\\
  &\longrightarrow R\Hom\bigl(R\Gamma_c(X,\Lambda),R\Gamma_c(X,\Lambda)\bigr),
 \end{align*}
 where the second arrow is induced from $\pr_{2*}\colon R\Gamma_c(X\times_SX,\Lambda(d)[2d])\longrightarrow 
 R\Gamma_c(X,\Lambda)$.
 Hence we may divide the diagram above into two parts:
 \[\small
 \xymatrix{
 R\Gamma\bigl(X,\Lambda(d)[2d]\bigr)\Lotimes R\Gamma_c(X,\Lambda)\ar[d]^-{\cong}
 \ar@<-2pt>[r]^-{(\dagger)}& R\Hom\bigl(R\Gamma_c(X,\Lambda),R\Gamma_c(X,\Lambda(d)[2d])\Lotimes R\Gamma_c(X,\Lambda)\bigr)\ar[dd]^-{\cong}\\
 R\Gamma\bigl(\overline{X}\times_S\overline{X},j_*\Lambda\Lboxtimes j_!\Lambda(d)[2d]\bigr)\ar[d]\\
 R\Gamma\bigl(\overline{X}\times_S\overline{X},(j\times 1)_*(1\times j)_!\Lambda(d)[2d]\bigr)\ar[r]& R\Hom\bigl(R\Gamma_c(X,\Lambda),
 R\Gamma_c(X\times_SX,\Lambda(d)[2d])\bigr)
 }
 \]
 and
 \[
 \xymatrix{%
 R\Gamma\bigl(X,\Lambda(d)[2d]\bigr)\Lotimes R\Gamma_c(X,\Lambda)\ar@<-2pt>[r]\ar[d]^-{(\dagger)}&
 R\Hom\bigl(R\Gamma_c(X,\Lambda),\Lambda\bigr)\Lotimes R\Gamma_c(X,\Lambda)\ar[d]\\
 R\Hom\bigl(R\Gamma_c(X,\Lambda),R\Gamma_c(X,\Lambda(d)[2d])\Lotimes R\Gamma_c(X,\Lambda)\bigr)\ar[d]
 \ar[r]^-{\Tr_X\otimes \id}& R\Hom\bigl(R\Gamma_c(X,\Lambda),R\Gamma_c(X,\Lambda)\bigr)\\
 R\Hom\bigl(R\Gamma_c(X,\Lambda),R\Gamma_c(X\times_SX,\Lambda(d)[2d])\bigr)\lefteqn{.}\ar[ur]}
 \]
 Here $(\dagger)$ is induced from 
 \[
 R\Gamma_c(X,\Lambda)\Lotimes R\Gamma(X,\Lambda(d)[2d])\Lotimes R\Gamma_c(X,\Lambda)\yrightarrow{\cup\otimes \id} R\Gamma_c(X,\Lambda(d)[2d])\Lotimes R\Gamma_c(X,\Lambda).
 \]
 Let us prove the commutativity of the former. By the adjointness, this is equivalent to the commutativity of
 the following:
 \[\small
 \xymatrix{
 R\Gamma_c(X,\Lambda)\Lotimes R\Gamma\bigl(X,\Lambda(d)[2d]\bigr)\Lotimes R\Gamma_c(X,\Lambda)\ar[d]^-{\id\otimes (\cup \circ (\pr_1^*\otimes \pr_2^*))}
 \ar@<-2pt>[r]^-{\cup\otimes\id}& R\Gamma_c(X,\Lambda(d)[2d])\Lotimes R\Gamma_c(X,\Lambda)\ar[dd]_-{\cong}^-{\cup\circ (\pr_1^*\otimes\pr_2^*)}\\
 R\Gamma_c(X,\Lambda)\Lotimes R\Gamma\bigl(\overline{X}\times_S\overline{X},j_*\Lambda \Lboxtimes j_!\Lambda(d)[2d]\bigr)\ar[d]
 \\
 R\Gamma_c(X,\Lambda)\Lotimes R\Gamma\bigl(\overline{X}\times_S\overline{X},(j\times 1)_*(1\times j)_!\Lambda(d)[2d]\bigr)\ar[r]^-{\cup\circ (\pr_1^*\otimes\id)}& 
 R\Gamma_c(X\times_SX,\Lambda(d)[2d])\lefteqn{.}
 }
 \]
 It easily follows from the associativity of cup products.

 Next we will prove the commutativity of the latter. The triangle is commutative,
 since the trace map for a smooth morphism is compatible with base change.
 By the adjointness, the commutativity of the rectangle is equivalent to that of the following:
 \[\small
 \xymatrix{%
 R\Gamma_c(X,\Lambda)\Lotimes R\Gamma\bigl(X,\Lambda(d)[2d]\bigr)\Lotimes R\Gamma_c(X,\Lambda)\ar@<-2pt>[r]\ar[d]^-{\cup\otimes \id}&
 R\Gamma_c(X,\Lambda)\Lotimes R\Hom\bigl(R\Gamma_c(X,\Lambda),\Lambda\bigr)\Lotimes R\Gamma_c(X,\Lambda)\ar[d]^-{\mathrm{ev}\otimes \id}\\
 R\Gamma_c\bigl(X,\Lambda(d)[2d]\bigr)\Lotimes R\Gamma_c(X,\Lambda)\ar[r]^-{\Tr_X\otimes \id}& R\Gamma_c(X,\Lambda)\lefteqn{.}
 }
 \]
 Since it is obtained from the diagram in Proposition \ref{prop:commutativity-part-1} by taking tensor products with 
 $R\Gamma_c(X,\Lambda)$, it is commutative.
\end{prf}

Now the proof of Proposition \ref{prop:LTF-unlocalized} is complete.

\subsection{Localization}\label{subsec:localization}
Let the notation be the same as in the previous subsection.
In this subsection, we will prove the main theorem in this paper, whose statement is the following:

\begin{thm}\label{thm:LTF-open-adic}
 Assume that for every $z\in \overline{\Gamma}\setminus \Gamma$, the points $\overline{\gamma}_1(z)$ and 
 $\overline{\gamma}_2(z)$ are distinct.
 Then we have
 \[
  \Tr\bigl(\gamma^*;R\Gamma_c(X,\Lambda)\bigr)=\#\Fix\gamma.
 \]
\end{thm}

\begin{rem}
 The condition in Theorem \ref{thm:LTF-open-adic} implies that 
 $\overline{\gamma}^{\,-1}(\Delta_{\overline{X}}\setminus \Delta_X)=\varnothing$,
 thus $\Fix\gamma=\Fix\overline{\gamma}$. In particular, $\Fix\gamma$ is proper over $S$
 and $\#\Fix \gamma$ makes sense.
 Note also that $\Delta_X\cap \gamma(\Gamma)=\Delta_{\overline{X}}\cap \overline{\gamma}(\overline{\Gamma})$,
 which will be used in the proof of Theorem \ref{thm:LTF-open-adic}.
\end{rem}

To prove Theorem \ref{thm:LTF-open-adic}, it suffices to compare the right hand side of
Proposition \ref{prop:LTF-unlocalized} with $\#\Fix \gamma$. 
The idea is to localize (or refine) the isomorphism
\[
 H^{2d}\bigl(\overline{X}\times_S\overline{X},j_*\Lambda \Lboxtimes j_!\Lambda(d)\bigr)\yrightarrow{\cong} H^{2d}\bigl(\overline{X}\times_S\overline{X},(j\times 1)_*(1\times j)_!\Lambda(d)\bigr)
\]
used in the previous subsection. 
Although the element $\cl(\gamma)$ lies in the local cohomology
$H^{2d}_{\overline{\gamma}(\overline{\Gamma})}(\overline{X}\times_S\overline{X},(j\times 1)_*(1\times j)_!\Lambda(d))$,
the map
\[
 H^{2d}_{\overline{\gamma}(\overline{\Gamma})}\bigl(\overline{X}\times_S\overline{X},j_*\Lambda \Lboxtimes j_!\Lambda(d)\bigr)\longrightarrow H^{2d}_{\overline{\gamma}(\overline{\Gamma})}\bigl(\overline{X}\times_S\overline{X},(j\times 1)_*(1\times j)_!\Lambda(d)\bigr)
\]
is not necessary an isomorphism (\cf Remark \ref{rem:isom-failure}). Thus we will slightly enlarge 
the closed subset $\overline{\gamma}(\overline{\Gamma})$ so that the morphism above becomes an isomorphism.
For the precise statement, see Proposition \ref{prop:localization}.

To do this, we need some preparation. First we show that $\overline{\gamma}_1(z)$ and 
$\overline{\gamma}_2(z)$ in Theorem \ref{thm:LTF-open-adic}
 can be separated by closed constructible subsets.

\begin{prop}\label{prop:separated-by-closed-constr}
 Under the assumption in Theorem \ref{thm:LTF-open-adic}, for every $z\in\overline{\Gamma}\setminus\Gamma$,
 the points $\overline{\gamma}_1(z)$ and $\overline{\gamma}_2(z)$ can be separated
 by closed constructible subsets of $\overline{X}$; namely, there exist closed constructible subsets
 $W_1$ and $W_2$ of $\overline{X}$ such that $\overline{\gamma}_1(z)\in W_1$, 
 $\overline{\gamma}_2(z)\in W_2$ and $W_1\cap W_2=\varnothing$. 
\end{prop}

\begin{rem}\label{rem:proper-spectral}
 The underlying topological space $\lvert \overline{X}\rvert$ of a quasi-compact quasi-separated 
 pseudo-adic space $\overline{X}$ is spectral.
 Indeed, it is by definition quasi-compact, quasi-separated and locally pro-constructible in the locally spectral space
 $(\overline{X})_{-}$ (\cf \cite[Remark 2.1 (iv)]{MR1207303}).
 Therefore, a subset $W$ of $\overline{X}$ is closed constructible if and only if $\overline{X}\setminus W$ is
 quasi-compact open (\cf \cite[Remark 2.1 (i)]{MR1207303}).
\end{rem}

\begin{rem}\label{rem:separated-by-closed-constr}
 For two points $x$, $y$ in a spectral space $X$, they can be separated by closed constructible subsets
 if and only if $\overline{\{x\}}\cap \overline{\{y\}}=\varnothing$. 
 Indeed, if $\overline{\{x\}}\cap \overline{\{y\}}=\varnothing$, we can cover $\overline{\{y\}}$ by
 finitely many quasi-compact open subsets $U_1,\ldots,U_m$ such that $U_i\cap \overline{\{x\}}=\varnothing$.
 Then, the complement $W_1$ of $U_1\cup \cdots \cup U_m$ is a closed constructible subset of $X$
 containing $x$. Since $W_1\cap \overline{\{y\}}=\varnothing$, we can cover $W_1$ by
 finitely many quasi-compact open subsets $V_1,\ldots,V_n$ such that $V_i\cap \overline{\{y\}}=\varnothing$.
 Let $W_2$ be the complement of $V_1\cup \cdots \cup V_n$. It is a closed constructible subset of $X$
 such that $y\in W_2$ and $W_1\cap W_2=\varnothing$.
 The converse direction is obvious.
\end{rem}

Since $\overline{\gamma}_1, \overline{\gamma}_2\colon \overline{\Gamma}\setminus \Gamma\longrightarrow \overline{X}$
are closed maps, Proposition \ref{prop:separated-by-closed-constr} is reduced to the following lemma:

\begin{lem}
 Let $X$, $Y$ be quasi-compact quasi-separated analytic pseudo-adic spaces and $f$, $g$ closed morphisms
 between them. Assume that for every $x\in X$, $f(x)\neq g(x)$.
 Then, for every $x\in X$, $f(x)$ and $g(x)$ can be separated by closed constructible subsets of $Y$.
\end{lem}

\begin{prf}
 Let $T$ be the subset of $X$ consisting of $x$ such that $f(x)$ and $g(x)$ can be separated by
 closed constructible subsets of $Y$. Then $T$ is an ind-constructible subset of $X$.
 Indeed, for each $x\in T$, we can find two closed constructible subsets $V_x$ and $W_x$ of $Y$
 separating $f(x)$ and $g(x)$, and then $f^{-1}(V_x)\cap g^{-1}(W_x)$ is contained in $T$.
 Therefore, $T$ is a union of constructible subsets $f^{-1}(V_x)\cap g^{-1}(W_x)$ for $x\in X$
 (recall that $f$ and $g$ are spectral \cite[Proposition 3.8 (iii), (iv)]{MR1207303}).

 Hence $X\setminus T$ is a pro-constructible subset of $X$,
 and thus it is a spectral space.
 Assume $X\setminus T$ is non-empty. Then it contains a closed point $z$
 (recall that every spectral space is homeomorphic to $\Spec A$ for some ring $A$
 \cite[Theorem 6]{MR0251026}). By Remark \ref{rem:separated-by-closed-constr}, we can find
 $y\in \overline{\{f(z)\}}\cap \overline{\{g(z)\}}$. By \cite[Lemma 1.1.10 i)]{MR1734903},
 we have either $f(z)\succ g(z)$ or $g(z)\succ f(z)$, where $\succ$ denotes the specialization relation
 as in \cite[p.~41]{MR1734903}. We may assume that $f(z)\succ g(z)$ without loss of generality.
 Since $f$ is closed, there exists $x\in X$ such that $z\succ x$ and $f(x)=g(z)$.
 Then we have $f(x)=g(z)\succ g(x)$, which implies $x\notin T$. 
 As $z$ is closed in $X\setminus T$ and
 $x$ is a specialization of $z$ inside $X\setminus T$, we have $x=z$, and thus $f(x)=g(x)$.
 This contradicts to our original assumption.

 Therefore we conclude that $X=T$, as desired.
\end{prf}

Next we consider the map between local cohomology induced by 
$\tau\colon j_*\Lambda \Lboxtimes j_!\Lambda\longrightarrow (j\times 1)_*(1\times j)_!\Lambda$.

\begin{lem}\label{lem:square-isom-open}
 Let $U$ and $V$ be quasi-compact open subsets of $\overline{X}$. Then the map
 $R\Gamma_{U\times_SV}(\overline{X}\times_S\overline{X},j_*\Lambda \Lboxtimes j_!\Lambda)\longrightarrow R\Gamma_{U\times_SV}(\overline{X}\times_S\overline{X},(j\times 1)_*(1\times j)_!\Lambda)$
 induced by $\tau$ is an isomorphism.
\end{lem} 

\begin{prf}
 By easy observation, we have
 \begin{align*}
  (j_*\Lambda \Lboxtimes j_!\Lambda)\vert_{U\times_SV}&\cong j'_*\Lambda \Lboxtimes j''_!\Lambda,\\
  \bigl((j\times 1)_*(1\times j)_!\Lambda\bigr)\bigr\vert_{U\times_SV}&\cong (j'\times 1)_*(1\times j'')_!\Lambda,
 \end{align*}
 where $j'$ (resp.\ $j''$) is the open immersion $U\cap X\hooklongrightarrow U$
 (resp.\ $V\cap X\hooklongrightarrow V$). Hence we have
 \begin{align*}
  R\Gamma_{U\times_SV}(\overline{X}\times_S\overline{X},j_*\Lambda \Lboxtimes j_!\Lambda)&=
  R\Gamma\bigl(U\times_SV,j'_*\Lambda \Lboxtimes j''_!\Lambda\bigr),\\
  R\Gamma_{U\times_SV}\bigl(\overline{X}\times_S\overline{X},(j\times 1)_*(1\times j)_!\Lambda\bigr)&=
  R\Gamma\bigl(U\times_SV,(j'\times 1)_*(1\times j'')_!\Lambda\bigr).
 \end{align*}

 Note that $U\cap X$ and $V\cap X$ are also quasi-compact, since $\overline{X}$ is quasi-separated.
 Therefore by Corollary \ref{cor:Kunneth-square-shape},
 the canonical map
 \[
 R\Gamma(U\times_SV,j'_*\Lambda \Lboxtimes j''_!\Lambda)\longrightarrow R\Gamma\bigl(U\times_SV,(j'\times 1)_*(1\times j'')_!\Lambda\bigr)
 \]
 is an isomorphism. Now the proof is complete.
\end{prf}

\begin{cor}\label{cor:square-isom-open-union}
 Let $U_1,\ldots,U_m, U'_1,\ldots,U'_m$ be quasi-compact open subsets of $X$ and put $V=\bigcup_{i=1}^mU_i\times_SU'_i$. 
 Then the map $R\Gamma_V(\overline{X}\times_S\overline{X},j'_*\Lambda \Lboxtimes j''_!\Lambda)\longrightarrow R\Gamma_V(\overline{X}\times_S\overline{X},(j\times 1)_*(1\times j)_!\Lambda)$ induced by $\tau$ is an isomorphism.
\end{cor}

\begin{prf}
 Recall that the underlying topological space of $U_i\times_SU'_i$ is equal to
 $\pr_1^{-1}(U_i)\cap \pr_2^{-1}(U'_i)$.
 Therefore, for every subset $I\subset \{1,\ldots,m\}$,
 we have $\bigcap_{i\in I}U_i\times_SU'_i=(\bigcap_{i\in I}U_i)\times_S(\bigcap_{i\in I}U'_i)$.
 Then an easy Mayer-Vietoris argument we may reduce to the case where $m=1$, which is already proven in 
 the previous lemma.
\end{prf}

\begin{cor}\label{cor:square-isom-closed}
 Let $V\subset \overline{X}\times_S\overline{X}$ be an open subset of the type in the above corollary. Put $W=(\overline{X}\times_S\overline{X})\setminus V$.
 Then the map $R\Gamma_W(\overline{X}\times_S\overline{X},j_*\Lambda \Lboxtimes j_!\Lambda)\longrightarrow R\Gamma_W(\overline{X}\times_S\overline{X},(j\times 1)_*(1\times j)_!\Lambda)$ induced by $\tau$ is an isomorphism.
\end{cor}

\begin{prf}
 Clear from the distinguished triangle 
 \[
  R\Gamma_V(\overline{X}\times_S\overline{X},-)\longrightarrow R\Gamma(\overline{X}\times_S\overline{X},-)\longrightarrow R\Gamma_W(\overline{X}\times_S\overline{X},-)\yrightarrow{+1}
 \]
 and Corollary \ref{cor:square-isom-open-union}.
\end{prf}

\begin{cor}\label{cor:square-isom-closed-union}
 Let $W_1,\ldots,W_m, W'_1,\ldots,W'_m$ be closed constructible subsets of $\overline{X}$ and put 
 $Z=\bigcup_{i=1}^mW_i\times_SW'_i$. 
 Then the map $R\Gamma_Z(\overline{X}\times_S\overline{X},j_*\Lambda \Lboxtimes j_!\Lambda)\longrightarrow R\Gamma_Z(\overline{X}\times_S\overline{X},(j\times 1)_*(1\times j)_!\Lambda)$ induced by $\tau$ is an isomorphism.
\end{cor}

\begin{prf}
 In the same way as in the proof of Corollary \ref{cor:square-isom-open-union}, we can reduce to the case 
 where $m=1$. This is the special case of Corollary \ref{cor:square-isom-closed},
 since $W_1\times_S W'_1$ is the complement of $(W_1^c\times_S \overline{X})\cup (\overline{X}\times_S W_1'^c)$
 (\cf Remark \ref{rem:proper-spectral}).
\end{prf}

Now we use the condition in Theorem \ref{thm:LTF-open-adic}. 
Proposition \ref{prop:separated-by-closed-constr} tells us that
for every $y\in \overline{\gamma}(\overline{\Gamma})\setminus (X\times_SX)$ there exist
closed constructible subsets $W_y$ and $W'_y$ of $\overline{X}$ such that
$y\in W_y\times_SW'_y$ and $W_y\cap W'_y=\varnothing$.
Since $\overline{\gamma}(\overline{\Gamma})\setminus (X\times_SX)$ is closed and $W_y\times_SW'_y$ is 
constructible in $\overline{X}\times_S\overline{X}$, we may choose finitely many
$W_1\times_SW_1',\ldots,W_m\times_SW_m'$ among $\{W_y\times_SW_y'\}$
so that they cover $\overline{\gamma}(\overline{\Gamma})\setminus (X\times_SX)$. 
Indeed, if we endow $\overline{X}\times_S\overline{X}$ with the patch topology,
$\overline{\gamma}(\overline{\Gamma})\setminus (X\times_SX)$ becomes compact and $W_y\times_SW'_y$ becomes open
(\cf Remark \ref{rem:proper-spectral}, \cite[\S 2]{MR0251026}). 

Put $Z=\bigcup_{i=1}^mW_i\times_SW_i'$.
The following proposition is crucial for the proof of Theorem \ref{thm:LTF-open-adic}:

\begin{prop}\label{prop:localization}
 The map
 \[
 R\Gamma_{\overline{\gamma}(\overline{\Gamma})\cup Z}(\overline{X}\times_S\overline{X},j_*\Lambda \Lboxtimes j_!\Lambda)\longrightarrow R\Gamma_{\overline{\gamma}(\overline{\Gamma})\cup Z}\bigl(\overline{X}\times_S\overline{X},(j\times 1)_*(1\times j)_!\Lambda\bigr)
 \]
 induced by $\tau$ is an isomorphism.
\end{prop}

\begin{prf}
 Put $V=(\overline{X}\times_S\overline{X})\setminus Z$.
 By the distinguished triangle 
 \[
  R\Gamma_Z(\overline{X}\times_S\overline{X},-)\longrightarrow R\Gamma_{\overline{\gamma}(\overline{\Gamma})\cup Z}(\overline{X}\times_S\overline{X},-)\longrightarrow R\Gamma_{\overline{\gamma}(\overline{\Gamma})\setminus Z}\bigl((V,(-)\vert_V\bigr)\yrightarrow{+1}
 \]
 and Corollary \ref{cor:square-isom-closed-union}, it suffices to show that
 the map
 \[
 R\Gamma_{\overline{\gamma}(\overline{\Gamma})\setminus Z}\bigl(V,(j_*\Lambda \Lboxtimes j_!\Lambda)\vert_V\bigr)\longrightarrow R\Gamma_{\overline{\gamma}(\overline{\Gamma})\setminus Z}\bigl(V,((j\times 1)_*(1\times j)_!\Lambda)\vert_V\bigr)
 \]
 is an isomorphism. By the assumption on $Z$, $\overline{\gamma}(\overline{\Gamma})\setminus Z$ has
 an open neighborhood $V\cap (X\times_SX)$ on which
 $\tau\colon j_*\Lambda \Lboxtimes j_!\Lambda\longrightarrow (j\times 1)_*(1\times j)_!\Lambda$ is an
 isomorphism. Thus the map above is also an isomorphism.
\end{prf}

 We define the element $[\gamma]_Z$ of $H^{2d}_{\overline{\gamma}(\overline{\Gamma})\cup Z}(\overline{X}\times_S\overline{X},j_*\Lambda \Lboxtimes j_!\Lambda(d))$ as
 the image of $\cl(\gamma)$ under the composite of following maps:
 \begin{align*}
  H^{2d}_{\gamma(\Gamma)}\bigl(X\times_SX,\Lambda(d)\bigr)&\cong
  H^{2d}_{\overline{\gamma}(\overline{\Gamma})}\bigl(\overline{X}\times_S\overline{X},(j\times 1)_*(1\times j)_!\Lambda(d)\bigr)\\
  &\longrightarrow H^{2d}_{\overline{\gamma}(\overline{\Gamma})\cup Z}\bigl(\overline{X}\times_S\overline{X},(j\times 1)_*(1\times j)_!\Lambda(d)\bigr)\\
  &\yleftarrow{\cong} H^{2d}_{\overline{\gamma}(\overline{\Gamma})\cup Z}\bigl(\overline{X}\times_S\overline{X},j_*\Lambda \Lboxtimes j_!\Lambda(d)\bigr).
 \end{align*}
 By the definition, the image of $[\gamma]_Z$ under the natural map
 \[
  H^{2d}_{\overline{\gamma}(\overline{\Gamma})\cup Z}\bigl(\overline{X}\times_S\overline{X},j_*\Lambda \Lboxtimes j_!\Lambda(d)\bigr)\longrightarrow H^{2d}\bigl(\overline{X}\times_S\overline{X},j_*\Lambda \Lboxtimes j_!\Lambda(d)\bigr)
 \]
 coincides with $[\gamma]$.

 Since $\Delta_{\overline{X}}\cap Z=\varnothing$, we have $\Delta_{\overline{X}}\cap (\overline{\gamma}(\overline{\Gamma})\cup Z)=\Delta_X\cap \gamma(\Gamma)$. Therefore $\overline{\delta}$ induces the maps
 \[
  H^{2d}_{\overline{\gamma}(\overline{\Gamma})\cup Z}\bigl(\overline{X}\times_S\overline{X},j_*\Lambda \Lboxtimes j_!\Lambda(d)\bigr)\yrightarrow{\overline{\delta}^*} H^{2d}_{\Delta_X\cap \gamma(\Gamma)}\bigl(\overline{X},j_!\Lambda(d)\bigr)\yrightarrow{\cong}
 H^{2d}_{\Delta_X\cap \gamma(\Gamma)}\bigl(X,\Lambda(d)\bigr).
 \]
 We denote the image of $[\gamma]_Z$ under the maps above by $\overline{\delta}^*([\gamma]_Z)$.

\begin{lem}\label{lem:class-local}
 The image of $\overline{\delta}^*([\gamma]_Z)$ under the canonical map
 $H^{2d}_{\Delta_X\cap\gamma(\Gamma)}(X,\Lambda(d))\longrightarrow H^{2d}_c(X,\Lambda(d))$
 coincides with $\overline{\delta}^*[\gamma]$.
\end{lem} 

\begin{prf}
 Clear from the following commutative diagram:
 \[
 \xymatrix{%
 H^{2d}_{\overline{\gamma}(\overline{\Gamma})\cup Z}\bigl(\overline{X}\times_S\overline{X},j_*\Lambda \Lboxtimes j_!\Lambda(d)\bigr)\ar[r]\ar[d]^-{\overline{\delta}^*}&
 H^{2d}\bigl(\overline{X}\times_S\overline{X},j_*\Lambda \Lboxtimes j_!\Lambda(d)\bigr)\ar[d]^-{\overline{\delta}^*}\\
 H^{2d}_{\Delta_X\cap\gamma(\Gamma)}\bigl(\overline{X},j_!\Lambda(d)\bigr)\ar[r]\ar[d]^-{\cong}&
 H^{2d}\bigl(\overline{X},j_!\Lambda(d)\bigr)\ar@{=}[d]\\
 H^{2d}_{\Delta_X\cap\gamma(\Gamma)}\bigl(X,\Lambda(d)\bigr)\ar[r]& H^{2d}_c\bigl(X,\Lambda(d)\bigr)\lefteqn{.}
 }
 \]
\end{prf}

\begin{lem}\label{lem:cl-gamma_Z}
 The image of $\cl(\gamma)$ under the map
 $\delta^*\colon H^{2d}_{\gamma(\Gamma)}(X\times_SX,\Lambda(d))\longrightarrow
 H^{2d}_{\Delta_X\cap \gamma(\Gamma)}(X,\Lambda(d))$
 coincides with $\overline{\delta}^*([\gamma]_Z)$.
\end{lem}

\begin{prf}
 Clear from the following commutative diagram:
 \[
 \scriptsize\xymatrix{%
 H^{2d}_{\overline{\gamma}(\overline{\Gamma})}\bigl(\overline{X}\times_S\overline{X},(j\times 1)_*(1\times j)_!\Lambda(d)\bigr)\ar[r]^-{\quad\vert_{X\times_SX}}_-{\cong}\ar[d]&H^{2d}_{\gamma(\Gamma)}\bigl(X\times_SX,\Lambda(d)\bigr)\ar[r]^-{\delta^*}\ar[d]&
 H^{2d}_{\Delta_X\cap \gamma(\Gamma)}\bigl(X,\Lambda(d)\bigr)\ar@{=}[d]\\
 H^{2d}_{\overline{\gamma}(\overline{\Gamma})\cup Z}\bigl(\overline{X}\times_S\overline{X},(j\times 1)_*(1\times j)_!\Lambda(d)\bigr)\ar[r]^-{\quad\vert_{X\times_SX}}&H^{2d}_{(X\times_SX)\cap (\overline{\gamma}(\overline{\Gamma})\cup Z)}\bigl(X\times_SX,\Lambda(d)\bigr)\ar[r]^-{\delta^*}&
 H^{2d}_{\Delta_X\cap \gamma(\Gamma)}\bigl(X,\Lambda(d)\bigr)\ar@{=}[d]\\
 H^{2d}_{\overline{\gamma}(\overline{\Gamma})\cup Z}\bigl(\overline{X}\times_S\overline{X},j_*\Lambda \Lboxtimes j_!\Lambda(d)\bigr)\ar[u]_{\cong}\ar[r]^-{\overline{\delta}^*}&
 H^{2d}_{\Delta_X\cap \gamma(\Gamma)}\bigl(\overline{X},j_!\Lambda(d)\bigr)\ar[r]^-{\quad\vert_X}_-{\cong}& H^{2d}_{\Delta_X\cap \gamma(\Gamma)}\bigl(X,\Lambda(d)\bigr)\lefteqn{.}
 }
 \]
\end{prf}

\begin{prf}[of Theorem \ref{thm:LTF-open-adic}]
 By Remark \ref{rem:Fix-via-cl}, $\#\Fix\gamma$ is the image of $\cl(\gamma)$ under the maps
 \[
  H^{2d}_{\gamma(\Gamma)}\bigl(X\times_SX,\Lambda(d)\bigr)\yrightarrow{\delta^*}
 H^{2d}_{\Delta_X\cap \gamma(\Gamma)}\bigl(X,\Lambda(d)\bigr)\longrightarrow 
 H^{2d}_c\bigl(X,\Lambda(d)\bigr)\yrightarrow{\Tr_X} \Lambda.
 \]
 By Lemma \ref{lem:class-local} and Lemma \ref{lem:cl-gamma_Z}, 
 it coincides with $\Tr_X(\overline{\delta}^*[\gamma])$.
 Therefore by Proposition \ref{prop:LTF-unlocalized}, we conclude that
 $\#\Fix \gamma=\Tr(\gamma^*;R\Gamma_c(X,\Lambda))$.
\end{prf}

\begin{rem}
 So far, we considered the case of torsion coefficient. However, at least when the characteristic of $k$ is $0$
 (\cf \cite[Theorem 3.1]{MR1626021}),
 we may obtain the Lefschetz trace formula for $\ell$-adic coefficient simply by taking projective limit.
 For the definition of $\#\Fix \gamma$ for the $\ell$-adic case, see Remark \ref{rem:Fix-l-adic}.
\end{rem}

\subsection{Lefschetz trace formula for open adic curves}\label{sec:LTF-curve}
In this subsection, we will establish a Lefschetz trace formula for quasi-compact smooth adic curve
by using the same idea as in the proof of Theorem \ref{thm:LTF-open-adic}.
Let $X$ be a $1$-dimensional quasi-compact adic space which is separated and smooth over $S$,
and $j\colon X\hooklongrightarrow X^c$ the universal compactification over $S$.
It is known that $\partial X:=X^c\setminus X$ is a finite discrete set (\cite[Lemma 5.12]{MR1856262}).
In particular, every $x\in \partial X$ is a closed constructible subset of $X^c$.

Let $f\colon X\longrightarrow X$ be a proper morphism over $S$ and $f^c\colon X^c\longrightarrow X^c$ the induced morphism.
Since $f$ is proper and $X$ is dense in $X^c$, we have $(f^c)^{-1}(\partial X)=\partial X$.
We will use the notation in Example \ref{exa:self-mor}.
Assume that $\Fix f$ is proper over $S$; thus $\#\Fix f$ can be defined.
Put $\partial X^{\fix}:=\{x\in \partial X\mid f^c(x)=x\}$. For $x\in \partial X^{\fix}$, we will define ``the contribution from $x$''
in the Lefschetz trace formula for $X$.

\begin{defn}
 Put $\Gamma'=\Gamma_f\cup \bigcup_{x\in \partial X}(f^c(x)\times_Sx)$.
 Since $x$ is a closed constructible subset of $X^c$, the natural map
 \[
  H^2_{\Gamma'}\bigl(X^c\times_SX^c,j_*\Lambda\Lotimes j_!\Lambda(1)\bigr)\longrightarrow H^2_{\Gamma'}\bigl(X^c\times_SX^c,(j\times 1)_*(1\times j)_!\Lambda(1)\bigr)
 \]
 is an isomorphism (\cf Proposition \ref{prop:localization}).
 We denote by $[f]_{\partial}$ the element of $H^2_{\Gamma'}(X^c\times_SX^c,j_*\Lambda\Lotimes j_!\Lambda(1))$ that is mapped to
 \[
  \cl(\gamma_f)\in H^2_{\Gamma_f}\bigl(X\times_SX,\Lambda(1)\bigr)=H^2_{\Gamma'}\bigl(X^c\times_SX^c,(j\times 1)_*(1\times j)_!\Lambda(1)\bigr).
 \]
 On the other hand, since $\Fix f$ is proper over $S$, it is a closed subset of $X^c$. Therefore we have
 $\Delta_{X^c}\cap \Gamma'=\Fix f\amalg \coprod_{x\in \partial X^{\fix}}x$ as a topological space, and thus
 $H^2_{\Delta_{X^c}\cap \Gamma'}(X^c,j_!\Lambda(1))=H^2_{\Fix f}(X^c,j_!\Lambda(1))\oplus \bigoplus_{x\in \partial X^{\fix}}H^2_x(X^c,j_!\Lambda(1))$.
 For $x\in \partial X^{\fix}$, we define $\loc(x)$ as the image of $[f]_{\partial}$ under the composite of
 \begin{align*}
  &H^2_{\Gamma'}\bigl(X^c\times_SX^c,j_*\Lambda\Lboxtimes j_!\Lambda(1)\bigr)\yrightarrow{(\delta^c)^*}H^2_{\Delta_{X^c}\cap \Gamma'}\bigl(X^c,j_!\Lambda(1)\bigr)\\\
  &\qquad=H^2_{\Fix f}\bigl(X^c,j_!\Lambda(1)\bigr)\oplus \bigoplus_{x\in \partial X^{\fix}}H^2_x\bigl(X^c,j_!\Lambda(1)\bigr)
  \longrightarrow H^2_x\bigl(X^c,j_!\Lambda(1)\bigr)\\
  &\qquad\longrightarrow H^2\bigl(X^c,j_!\Lambda(1)\bigr)=H^2_c\bigl(X,\Lambda(1)\bigr)\yrightarrow{\Tr_X}\Lambda,
 \end{align*}
 where $\delta^c$ denotes the diagonal morphism for $X^c$.
\end{defn}

The following is our Lefschetz trace formula for adic curves:

\begin{thm}
 In the setting above, we have
 \[
  \Tr\bigl(f^*;R\Gamma_c(X,\Lambda)\bigr)=\#\Fix f+\sum_{x\in\partial X^{\fix}}\loc(x).
 \]
\end{thm}

\begin{prf}
 It is easy to see that the image of $[f]_{\partial}$ under the composite of
 \begin{align*}
  &H^2_{\Gamma'}\bigl(X^c\times_SX^c,j_*\Lambda\Lboxtimes j_!\Lambda(1)\bigr)\yrightarrow{(\delta^c)^*}H^2_{\Delta_{X^c}\cap \Gamma'}\bigl(X^c,j_!\Lambda(1)\bigr)\\
  &\qquad=H^2_{\Fix f}\bigl(X^c,j_!\Lambda(1)\bigr)\oplus \bigoplus_{x\in \partial X^{\fix}}H^2_x\bigl(X^c,j_!\Lambda(1)\bigr)
  \longrightarrow H^2_{\Fix f}\bigl(X^c,j_!\Lambda(1)\bigr)\\
  &\qquad\longrightarrow H^2\bigl(X^c,j_!\Lambda(1)\bigr)=H^2_c\bigl(X,\Lambda(1)\bigr)\yrightarrow{\Tr_X}\Lambda
 \end{align*}
 is equal to $\#\Fix f$ (\cf the proof of Lemma \ref{lem:cl-gamma_Z}). Therefore, the theorem immediately follows from
 Proposition \ref{prop:LTF-unlocalized}.
\end{prf}

\begin{rem}
 The formula in Theorem \ref{sec:LTF-curve} is very similar to Huber's trace formula for open curves 
 (\cite[Theorem 6.3]{MR1856262}). However, the definition of his local term is different from ours;
 that is given by purely algebraic manner and depends only on the homomorphism 
 induced on the valuation ring corresponding to $x\in \partial X^{\fix}$.
 The author expects that these two local terms coincide, and will consider this problem in his future work.
\end{rem}

\section{Lefschetz trace formula for formal schemes}\label{sec:LTF-formal}
In this section, we deduce a Lefschetz trace formula for formal schemes
from Theorem \ref{thm:LTF-open-adic}. For formal schemes, we will use the same notation as
in \cite[\S 4]{formalnearby}. Let us recall some of them.

Let $R$ be a complete discrete valuation ring with separably closed residue field
and $k$ an algebraic closure of the fraction field $F$ of $R$.
Put $\mathcal{S}=\Spf R$.
Let $\mathcal{X}$ be a quasi-compact special formal scheme which is separated over $\mathcal{S}$.
Then we can associate $\mathcal{X}$ with the adic spaces $t(\mathcal{X})_a$, $t(\mathcal{X})_{\eta}$
and $t(\mathcal{X})_{\overline{\eta}}$. The adic space $t(\mathcal{X})_a$ is an open adic subspace
of $t(\mathcal{X})$ consisting of analytic points of $t(\mathcal{X})$. It is quasi-compact.
The adic space $t(\mathcal{X})_{\eta}$ is the ``generic fiber'' of $t(\mathcal{X})$,
which is locally of finite type, separated and taut over $\Spa(F,R)$.
The adic space $t(\mathcal{X})_{\overline{\eta}}$ is the base change of $t(\mathcal{X})_{\eta}$
from $\Spa(F,R)$ to $S=\Spa(k,k^+)$. Note that $t(\mathcal{X})_{\eta}$ and
$t(\mathcal{X})_{\overline{\eta}}$ are not necessarily quasi-compact.
In the sequel, we write $X$, $X_\eta$ and $X_{\overline{\eta}}$ for $t(\mathcal{X})_a$, 
$t(\mathcal{X})_{\eta}$ and $t(\mathcal{X})_{\overline{\eta}}$, respectively.
On the other hand, we denote the special fiber of $\mathcal{X}$ (resp.\ $X$) by $\mathcal{X}_s$ (resp.\ $X_s$).

Let $\mathcal{T}$ be a finite set equipped with a partial order
and $\{\mathcal{Y}_\alpha\}_{\alpha\in\mathcal{T}}$ a family of closed formal subschemes
of $\mathcal{X}_s$ indexed by $\mathcal{T}$. We put 
$Y_\alpha=t(\mathcal{Y}_{\alpha})_a=t(\mathcal{Y}_\alpha)\times_{t(\mathcal{X})}X$.
We assume the following:

\begin{assump}\label{assump:formal-sch}
 \begin{enumerate}
  \item $X_s=\bigcup_{\alpha\in\mathcal{T}}Y_{\alpha}$.
  \item For $\alpha\in\mathcal{T}$, put $Y_{(\alpha)}=Y_\alpha\setminus\bigcup_{\beta>\alpha}Y_\beta$. Then, for $\alpha,\beta\in\mathcal{T}$ with $\alpha\neq \beta$,
	$Y_{(\alpha)}\cap Y_{(\beta)}=\varnothing$.
 \end{enumerate}
\end{assump}

\begin{exa}\label{exa:assump-scheme}
 Let $\mathsf{X}$ be a scheme which is separated of finite type over $\Spec R$ and
 $\{\mathsf{Y}_\alpha\}_{\alpha\in\mathcal{T}}$ a family of closed subschemes of the special fiber
 $\mathsf{X}_s$ of $\mathsf{X}$. Assume the following conditions:
 \begin{enumerate}
  \item $\mathsf{X}_s=\bigcup_{\alpha\in\mathcal{T}}\mathsf{Y}_{\alpha}$.
  \item For $\alpha\in\mathcal{T}$, put $\mathsf{Y}_{(\alpha)}=\mathsf{Y}_\alpha\setminus\bigcup_{\beta>\alpha}\mathsf{Y}_\beta$. Then, for $\alpha,\beta\in\mathcal{T}$ with $\alpha\neq \beta$,
	$\mathsf{Y}_{(\alpha)}\cap \mathsf{Y}_{(\beta)}=\varnothing$.
  \item There exists the unique maximal element $\alpha_0$ in $\mathcal{T}$.
 \end{enumerate}
 Denote the completion of $\mathsf{X}$ along $\mathsf{Y}_{\alpha_0}$ by $\mathcal{X}$ and
 put $\mathcal{Y}_{\alpha}=\mathsf{Y}_{\alpha}\times_{\mathsf{X}}\mathcal{X}$ for each $\alpha\in\mathcal{T}$.
 Then $\mathcal{X}$ and $\{\mathcal{Y}_{\alpha}\}_{\alpha\in\mathcal{T}\setminus\{\alpha_0\}}$
 satisfy the assumption above.
 Indeed, we have the natural morphism of locally ringed spaces
 $(X,\mathcal{O}_X)\longrightarrow \mathsf{X}\setminus\mathsf{Y}_{\alpha_0}$
 (\cf \cite[Remark 4.6 (iv)]{MR1306024}) such that
 the inverse image of $\mathsf{Y}_{\alpha}\setminus \mathsf{Y}_{\alpha_0}$ is equal to $Y_{\alpha}$.
\end{exa}

Let us consider an isomorphism $f\colon \mathcal{X}\yrightarrow{\cong}\mathcal{X}$ over $\mathcal{S}$.
We also denote the induced isomorphism $X\yrightarrow{\cong}X$ by the same symbol $f$.
The induced isomorphisms $X_{\eta}\yrightarrow{\cong}X_{\eta}$ and 
$X_{\overline{\eta}}\yrightarrow{\cong}X_{\overline{\eta}}$ are denoted by $f_{\eta}$ and $f_{\overline{\eta}}$,
respectively.
We will make the following assumption on $f$:

\begin{assump}\label{assump:morphism}
 There exists an order-preserving bijection $f\colon \mathcal{T}\yrightarrow{\cong}\mathcal{T}$
 and a system of closed constructible subsets $\{Y_\alpha(n)\}_{n\ge 1}$
 of $X$ for each $\alpha\in \mathcal{T}$ satisfying the following:
 \begin{enumerate}
  \item $Y_\alpha(n+1)\subset Y_\alpha(n)$ for every $n\ge 1$. 
  \item $\bigcap_{n\ge 1}Y_\alpha(n)=Y_\alpha$.
  \item $f(Y_\alpha(n))=Y_{f(\alpha)}(n)$ for every $\alpha\in\mathcal{T}$ and $n\ge 1$.
  \item $f(\alpha)\neq \alpha$ for every $\alpha\in\mathcal{T}$.
 \end{enumerate}
\end{assump}

\begin{rem}
 Later we will give some conditions for the existence of $\{Y_\alpha(n)\}_{n\ge 1}$
 in Assumption \ref{assump:morphism}.
 In fact, one of the following is sufficient (Proposition \ref{prop:suff-cond-isom},
 Proposition \ref{prop:suff-cond-top-fin-order}):
 \begin{itemize}
  \item For every $\alpha\in \mathcal{T}$, the isomorphism $f\colon \mathcal{X}\yrightarrow{\cong}\mathcal{X}$
	induces an isomorphism of formal schemes $\mathcal{Y}_\alpha\yrightarrow{\cong} \mathcal{Y}_{f(\alpha)}$.
  \item For every $\alpha\in \mathcal{T}$, the isomorphism $f\colon \mathcal{X}\yrightarrow{\cong}\mathcal{X}$
	induces a set-theoretic bijection $Y_\alpha\yrightarrow{\cong}Y_{f(\alpha)}$.
	Furthermore, for every ideal of definition $\mathcal{I}$ of $\mathcal{X}$, there exists an integer $N\ge 1$
	such that $f^N\equiv \id \pmod{\mathcal{I}}$.
 \end{itemize}
\end{rem}

Now we can state our Lefschetz trace formula for $\mathcal{X}$:

\begin{thm}\label{thm:LTF-formal}
 In addition to Assumption \ref{assump:formal-sch} and Assumption \ref{assump:morphism}, assume that
 $\mathcal{X}$ is locally algebraizable (\cite[Definition 3.18]{formalnearby}) and
 $X_{\eta}$ is partially proper and smooth over $\Spa(F,R)$. 
 Then $R\Gamma_c(X_{\overline{\eta}},\Lambda)$ is a perfect $\Lambda$-complex,
 $\Fix f_{\overline{\eta}}$ (\cf Example \ref{exa:self-mor}) is proper over $S$ and we have
 \[
  \Tr\bigl(f_{\overline{\eta}}^*;R\Gamma_c(X_{\overline{\eta}},\Lambda)\bigr)=\#\Fix f_{\overline{\eta}}.
 \]
\end{thm}

In order to prove this theorem, we need some preparation.
First we observe the finiteness of the cohomology of $X_{\overline{\eta}}$.

\begin{prop}\label{prop:nearby-finiteness}
 Let $\mathcal{X}$ be a quasi-compact special formal scheme which is locally algebraizable and
 separated over $\mathcal{S}$.
 Then $H^i_c(X_{\overline{\eta}},\Lambda)$ is a finitely generated $\Lambda$-module for every $i$.
 Thus $R\Gamma_c(X_{\overline{\eta}},\Lambda)$ is a perfect $\Lambda$-complex by Corollary \ref{cor:perfect-cpx}.
 More generally, let $L$ be a locally closed constructible subset of $X$.
 Then $H^i_c(L_{\overline{\eta}},\Lambda)$ is a finitely generated $\Lambda$-module for every $i$,
 and $R\Gamma_c(L_{\overline{\eta}},\Lambda)$ is a perfect $\Lambda$-complex.
\end{prop}

\begin{prf}
 We way assume that $\mathcal{X}$ is algebraizable.
 First we consider the case where $L$ is a quasi-compact open subset of $X$.
 Then there exists an admissible blow-up $\mathcal{X}'\longrightarrow \mathcal{X}$ and an open formal subscheme
 $\mathcal{U}'\subset \mathcal{X}'$ such that $L=t(\mathcal{U}')_a$. Since $\mathcal{U}'$ is algebraizable
 (\cf \cite[Lemma 7.1.4]{MR2309992}), replacing $\mathcal{X}$ by $\mathcal{U}'$, we may assume that $L=X$.
 Moreover we may assume that $\mathcal{X}$ is affine, and thus pseudo-compactifiable 
 (\cf \cite[Definition 4.21 i), Example 4.22 i)]{formalnearby}). Therefore, by \cite[Theorem 4.32]{formalnearby},
 we have $H^i_c(X_{\overline{\eta}},\Lambda)\cong H_c^i(\mathcal{X}_{\mathrm{red}},R\Psi_{\!\mathcal{X},c}\Lambda)$.
 On the other hand, by \cite[Proposition 3.20]{formalnearby}, $R\Psi_{\!\mathcal{X},c}\Lambda$ is constructible,
 and thus $H_c^i(\mathcal{X}_{\mathrm{red}},R\Psi_{\!\mathcal{X},c}\Lambda)$ is a finitely generated $\Lambda$-module.
 Hence $H^i_c(X_{\overline{\eta}},\Lambda)$ is also finitely generated, which concludes the proof for a quasi-compact
 open $L$.

 A general $L$ can be expressed as $L_2\setminus L_1$, where $L_1$ and $L_2$ are quasi-compact open subsets of $X$
 with $L_1\subset L_2$.
 Thus the proposition follows from the exact sequence
 \[
  H^i_c(L_{2,\overline{\eta}},\Lambda)\longrightarrow 
 H^i_c\bigl((L_2\setminus L_1)_{\overline{\eta}},\Lambda\bigr)
 \longrightarrow H^{i+1}_c(L_{1,\overline{\eta}},\Lambda).
 \]
\end{prf}

Now we use the notation in Assumption \ref{assump:morphism}.

\begin{lem}\label{lem:patch-nbhd}
 For every $\alpha\in\mathcal{T}$, $\{Y_\alpha(n)\}_{n\ge 1}$
 form a fundamental system of open neighborhoods of $Y_\alpha$ with respect to the patch topology of $X$.
\end{lem}

\begin{prf}
 Let $U$ be a subset of $X$ containing $Y_\alpha$, which is open in the patch topology.
 We will find $n\ge 1$ such that $Y_\alpha(n)\subset U$. 
 By Assumption \ref{assump:morphism} ii), $X\setminus U$ is covered by $\{X\setminus Y_\alpha(n)\}_{n\ge 1}$.
 Since $X\setminus U$ is compact and $X\setminus Y_\alpha(n)$ is an open subset of $X$
 with respect to the patch topology, there exists an integer $n\ge 1$ such that
 $X\setminus U\subset X\setminus Y_\alpha(n)$. In other words, $Y_\alpha(n)$ is contained in $U$.
\end{prf}

The following construction is crucial for the proof of Theorem \ref{thm:LTF-formal}:

\begin{lem}\label{lem:topology}
 We can find an integer $n_\alpha\ge 1$ for each $\alpha\in \mathcal{T}$ satisfying the following conditions:
 \begin{itemize}
  \item For every $\alpha\in\mathcal{T}$, $n_\alpha=n_{f(\alpha)}$.
  \item For $\alpha\in\mathcal{T}$, put $U_\alpha=Y_\alpha(n_\alpha)\setminus \bigcup_{\beta>\alpha}Y_\beta(n_\beta)$.
	Then we have $U_\alpha\cap U_\beta=\varnothing$ for every $\alpha,\beta\in\mathcal{T}$ with $\alpha\neq \beta$.
 \end{itemize}
\end{lem}

By the assumption, $Y(n)$ is an open subset of $X$ with respect to 
the patch topology. Therefore, the lemma is reduced to the following:

\begin{lem}
 Let $X$ be a compact topological space. Let $\mathcal{T}$ be a finite set equipped with a partial order
 and $\{Y_\alpha\}_{\alpha_\in\mathcal{T}}$ a family of closed subsets of $X$ indexed by $\mathcal{T}$.
 Put $Y_{(\alpha)}=Y_\alpha\setminus \bigcup_{\beta>\alpha}Y_\beta$ and assume that 
 $Y_{(\alpha)}\cap Y_{(\beta)}=\varnothing$ if $\alpha\neq \beta$.
 Let $f\colon \mathcal{T}\yrightarrow{\cong}\mathcal{T}$ be an order-preserving bijection.
 Assume that we are given a fundamental system of open neighborhoods $\{Y_\alpha(n)\}_{n\ge 1}$ of $Y_\alpha$
 for each $\alpha\in\mathcal{T}$ such that $Y_\alpha(n+1)\subset Y_\alpha(n)$ for every $n\ge 1$.
 Then, for every integer $N\ge 1$ we can find an integer $n_\alpha\ge N$ for each $\alpha\in \mathcal{T}$ 
 satisfying the following conditions:
 \begin{itemize}
  \item For every $\alpha\in\mathcal{T}$, $n_\alpha=n_{f(\alpha)}$.
  \item For $\alpha\in\mathcal{T}$, put $U_\alpha=Y_\alpha(n_\alpha)\setminus \bigcup_{\beta>\alpha}Y_\beta(n_\beta)$.
	Then we have $U_\alpha\cap U_\beta=\varnothing$ for every $\alpha,\beta\in\mathcal{T}$ with $\alpha\neq \beta$.
 \end{itemize}
\end{lem}

\begin{prf}
 Use the induction on the cardinality of $\mathcal{T}$. If $\mathcal{T}$ is empty, then the lemma is clear.
 Assume that $\mathcal{T}$ is non-empty.
 Take a maximal element $\alpha_0$ of $\mathcal{T}$ and put $\mathcal{T}_0=\{f^m(\alpha_0)\mid m\in\Z\}$.
 Note that every element of $\mathcal{T}_0$ is maximal.
 
 Consider an element $(\alpha,\beta)$ of $\mathcal{T}_0\times\mathcal{T}$ such that $\beta\nleq\alpha$.
 As $Y_\beta$ is contained in $\coprod_{\gamma\ge\beta}Y_{(\gamma)}$,
 $Y_\alpha=Y_{(\alpha)}$ does not intersect $Y_\beta$ by the assumption.
 Since $X$ is compact, two closed subsets $Y_\alpha$ and $Y_\beta$ can be separated by open neighborhoods.
 Thus we can find an integer $n\ge N$ such that $Y_\alpha(n)\cap Y_\beta(n)=\varnothing$.
 Since there are only finitely many such elements $(\alpha,\beta)$, we can take $n\ge N$
 such that $Y_\alpha(n)\cap Y_\beta(n)=\varnothing$ for every $\alpha\in\mathcal{T}_0$ and $\beta\in\mathcal{T}$
 with $\beta\nleq \alpha$. Put $n_\alpha=n$ for every $\alpha\in\mathcal{T}_0$.

 Put $X'=X\setminus \bigcup_{\alpha\in\mathcal{T}_0}Y_\alpha(n_\alpha)$, 
 $\mathcal{T}'=\mathcal{T}\setminus \mathcal{T}_0$, $Y'_\beta=X'\cap Y_\beta$ and $Y'_\beta(n)=X'\cap Y_\beta(n)$
 for $\beta\in\mathcal{T}'$. Note that $f$ induces an order-preserving bijection 
 $\mathcal{T}'\yrightarrow{\cong}\mathcal{T}'$, which is also denoted by $f$.
 Then, these satisfy the assumptions in the lemma. Therefore, by the induction hypothesis,
 we can find $n_\alpha\ge n_{\alpha_0}$ for each $\alpha\in\mathcal{T}'$.
 We will observe that $\{n_\alpha\}_{\alpha\in\mathcal{T}}$ satisfies the conditions in the lemma.
 The first condition $n_\alpha=n_{f(\alpha)}$ is clear from the construction.
 For the second condition, put $U'_\alpha=Y'_\alpha(n_\alpha)\setminus \bigcup_{\beta\in\mathcal{T'},\beta>\alpha}Y'_\beta(n_\beta)$ for $\alpha\in\mathcal{T}'$. Then,
 \begin{align*}
  U'_\alpha
  &=\Bigl(Y_\alpha(n_\alpha)\setminus \bigcup_{\beta\in\mathcal{T'},\beta>\alpha}Y_\beta(n_\beta)\Bigr)
 \setminus \bigcup_{\gamma\in\mathcal{T}_0}Y_\gamma(n_{\gamma})\\
  &=\Bigl(Y_\alpha(n_\alpha)\setminus \bigcup_{\beta\in\mathcal{T}',\beta>\alpha}Y_\beta(n_\beta)\Bigr)\setminus \bigcup_{\gamma\in\mathcal{T}_0,\gamma\ge \alpha}Y_\gamma(n_{\gamma})\\
  &=Y_\alpha(n_\alpha)\setminus \bigcup_{\beta\in\mathcal{T},\beta>\alpha}Y_\beta(n_\beta)=U_\alpha.
 \end{align*}
 The second equality follows from 
 $Y_\alpha(n_\alpha)\cap Y_\gamma(n_\gamma)\subset Y_\alpha(n_{\alpha_0})\cap Y_\gamma(n_{\alpha_0})=\varnothing$
 for $\gamma\in \mathcal{T}_0$ with $\gamma\ngeq \alpha$. 
 Let us take $\alpha, \beta\in \mathcal{T}$ with $\alpha\neq \beta$ and prove $U_\alpha\cap U_\beta=\varnothing$.
 If $\alpha,\beta\in\mathcal{T}_0$, $U_\alpha\cap U_\beta=Y_{\alpha}(n_{\alpha_0})\cap Y_{\beta}(n_{\alpha_0})=\varnothing$
 since $\beta\ngeq \alpha$. If $\alpha\in \mathcal{T}_0$ and $\beta\in \mathcal{T}'$, then
 $U_\alpha\cap U_\beta=Y_\alpha(n_\alpha)\cap U'_\beta\subset Y_\alpha(n_\alpha)\cap X'=\varnothing$.
 The case where $\alpha\in \mathcal{T}'$ and $\beta\in \mathcal{T}_0$ is similar.
 Finally if $\alpha,\beta\in \mathcal{T}'$, then $U_\alpha\cap U_\beta=U'_\alpha\cap U'_\beta=\varnothing$
 by the induction hypothesis.
 This completes the proof.
\end{prf}

Fix $\{n_\alpha\}_{\alpha\in\mathcal{T}}$ as in Lemma \ref{lem:topology} and
put $W=\bigcup_{\alpha\in\mathcal{T}}Y_\alpha(n_\alpha)$, $X_0=X\setminus W$. 
By Assumption \ref{assump:morphism} iii), we have $f(W)=W$ and $f(X_0)=X_0$.

\begin{prop}\label{prop:trace-identity}
 We have $\Tr(f_{\overline{\eta}}^*;R\Gamma_c(X_{\overline{\eta}},\Lambda))=\Tr(f^*_{\overline{\eta}};R\Gamma_c(X_{0,\overline{\eta}},\Lambda))$.
\end{prop}

\begin{lem}\label{lem:Tr-additivity}
 Let $L$, $L'$ be locally closed constructible subsets of $X$ such that $f(L)=L$, $f(L')=L'$
 and $L'$ is an open subset of $L$. Put $L''=L\setminus L'$.
 Then we have 
 \[
  \Tr\bigl(f_{\overline{\eta}}^*;R\Gamma_c(L_{\overline{\eta}},\Lambda)\bigr)
 =\Tr\bigl(f_{\overline{\eta}}^*;R\Gamma_c(L'_{\overline{\eta}},\Lambda)\bigr)
 +\Tr\bigl(f_{\overline{\eta}}^*;R\Gamma_c(L''_{\overline{\eta}},\Lambda)\bigr).
 \]
\end{lem}

\begin{prf}
 First note that $R\Gamma_c(L_{\overline{\eta}},\Lambda)$, $R\Gamma_c(L'_{\overline{\eta}},\Lambda)$ and
 $R\Gamma_c(L''_{\overline{\eta}},\Lambda)$ are perfect $\Lambda$-complexes by Proposition \ref{prop:nearby-finiteness},
 and the traces make sense.

 Let $j\colon (X_{\overline{\eta}},L'_{\overline{\eta}})\hooklongrightarrow (X_{\overline{\eta}},L_{\overline{\eta}})$
 and $i\colon (X_{\overline{\eta}},L''_{\overline{\eta}})\hooklongrightarrow (X_{\overline{\eta}},L_{\overline{\eta}})$
 be the natural immersions of pseudo-adic spaces. Consider the filtered sheaf 
 $\mathcal{F}=(F_1\mathcal{F}=j_!\Lambda\subset \Lambda=F_0\mathcal{F})$.
 We have a natural morphism $f_{\overline{\eta}}^*\mathcal{F}\longrightarrow \mathcal{F}$ of filtered sheaves,
 which induces a morphism
 $f_{\overline{\eta}}^*\colon R\Gamma_c(L_{\overline{\eta}},\mathcal{F})\longrightarrow R\Gamma_c(L_{\overline{\eta}},f_{\overline{\eta}}^*\mathcal{F})
 \longrightarrow R\Gamma_c(L_{\overline{\eta}},\mathcal{F})$
 in the filtered derived category of $\Lambda$-modules (\cf \cite[Chapitre V]{MR0491680}). 
 It is easy to see that the morphism on $\mathrm{gr}^0$ (resp.\ $\mathrm{gr}^1$) induced by $f_{\overline{\eta}}^*$
 coincides with the pull-back map $f_{\overline{\eta}}^*$ on 
 $R\Gamma_c(L_{\overline{\eta}},i_*\Lambda)=R\Gamma_c(L''_{\overline{\eta}},\Lambda)$
 (resp.\ $R\Gamma_c(L_{\overline{\eta}},j_!\Lambda)=R\Gamma_c(L'_{\overline{\eta}},\Lambda)$).
 Therefore the equality follows from \cite[Corollaire 3.7.7, Remarque 3.7.7.1]{MR0491680}.
\end{prf}

\begin{prf}[of Proposition \ref{prop:trace-identity}]
 By Lemma \ref{lem:Tr-additivity}, it suffices to show 
 $\Tr(f_{\overline{\eta}}^*;R\Gamma_c(W_{\overline{\eta}},\Lambda))=0$.
 Take a maximal element $\alpha_0$ of $\mathcal{T}$ and put $\mathcal{T}_0=\{f^m(\alpha_0)\mid m\in\Z\}$,
 $W_0=\bigcup_{\alpha\in\mathcal{T}_0}Y_\alpha(n_\alpha)=\coprod_{\alpha\in\mathcal{T}_0}U_\alpha$.
 Obviously $R\Gamma_c(W_{0,\overline{\eta}},\Lambda)=\bigoplus_{\alpha\in\mathcal{T}_0}R\Gamma_c(U_{\alpha,\overline{\eta}},\Lambda)$ and $f(W_0)=W_0$.
 Since $f(U_\alpha)=U_{f(\alpha)}$ and $f(\alpha)\neq \alpha$ by Assumption \ref{assump:morphism} iv),
 it is immediate to see $\Tr(f_{\overline{\eta}}^*;R\Gamma_c(W_{0,\overline{\eta}},\Lambda))=0$.

 Put $W'=W\setminus W_0$ and $\mathcal{T}'=\mathcal{T}\setminus \mathcal{T}_0$.
 Then $\Tr(f_{\overline{\eta}}^*;R\Gamma_c(W_{\overline{\eta}},\Lambda))=\Tr(f_{\overline{\eta}}^*;R\Gamma_c(W'_{\overline{\eta}},\Lambda))$ by Lemma \ref{lem:Tr-additivity}.
 If $\mathcal{T}'$ is non-empty, take a maximal element $\alpha_1$ of $\mathcal{T}'$
 and put $\mathcal{T}_1=\{f^m(\alpha_1)\mid m\in\Z\}$,
 $W_1=\bigcup_{\alpha\in\mathcal{T}_1}Y_\alpha(n_\alpha)\setminus W_0=\coprod_{\alpha\in\mathcal{T}_1}U_\alpha$.
 In the same way as above, we can prove that 
 $\Tr(f_{\overline{\eta}}^*;R\Gamma_c(W_{1,\overline{\eta}},\Lambda))=0$.
 Put $W''=W'\setminus W_1$. 
 Then $\Tr(f_{\overline{\eta}}^*;R\Gamma_c(W'_{\overline{\eta}},\Lambda))=\Tr(f_{\overline{\eta}}^*;R\Gamma_c(W''_{\overline{\eta}},\Lambda))$ by Lemma \ref{lem:Tr-additivity} (note that $W_1$ is closed in $W'$).
 We repeat this procedure to obtain $\Tr(f_{\overline{\eta}}^*;R\Gamma_c(W_{\overline{\eta}},\Lambda))=0$.
\end{prf}

Next lemma ensures that we may apply Theorem \ref{thm:LTF-open-adic} to $X_{0,\overline{\eta}}$.

\begin{lem}\label{lem:verify-assump}
 \begin{enumerate}
  \item The adic space $X_0$ is a quasi-compact open adic subspace of $X_{\eta}$.
	In particular, $X_{0,\overline{\eta}}$ is smooth, separated of finite type over $S$.
  \item For $x\in X_{\overline{\eta}}\setminus X_{0,\overline{\eta}}$, we have $x\neq f_{\overline{\eta}}(x)$.
  \item The closed adic subspace $\Fix f_{\overline{\eta}}$ of $X_{\overline{\eta}}$ is contained 
	in $X_{0,\overline{\eta}}$.
	In particular, $\Fix f_{\overline{\eta}}=\Fix (f_{\overline{\eta}}\vert_{X_{0,\overline{\eta}}})$ is proper over $S$ and
	$\#\Fix f_{\overline{\eta}}=\#\Fix (f_{\overline{\eta}}\vert_{X_{0,\overline{\eta}}})$.
 \end{enumerate}
\end{lem}

\begin{prf}
 i) Since $X$ is a spectral space and $W$ is a closed constructible subset of $X$,
 $X_0$ is a quasi-compact open subset of $X$.
 On the other hand, $X_0\subset X\setminus \bigcup_{\alpha\in\mathcal{T}}Y_\alpha=X_\eta$ by 
 Assumption \ref{assump:formal-sch} i). Thus $X_0$ is a quasi-compact open subset of $X_\eta$.

 \noindent ii) As $x\in W_{\overline{\eta}}=\coprod_{\alpha\in\mathcal{T}}U_{\alpha,\overline{\eta}}$,
 there exists $\alpha\in\mathcal{T}$ such that $x\in U_{\alpha,\overline{\eta}}$. 
 By Assumption \ref{assump:morphism} iii), we have $f_{\overline{\eta}}(x)\in U_{f(\alpha),\overline{\eta}}$.
 Since $U_{\alpha,\overline{\eta}}\cap U_{f(\alpha),\overline{\eta}}=\varnothing$
 due to Assumption \ref{assump:morphism} iv),
 $x$ and $f_{\overline{\eta}}(x)$ are distinct.

 \noindent iii) Clear from ii) and Proposition \ref{prop:Fix-etale}.
\end{prf}

\begin{prf}[of Theorem \ref{thm:LTF-formal}]
 Since $X_{\overline{\eta}}$ is partially proper and taut over $S$ by the assumption,
 the closure $\overline{X_{0,\overline{\eta}}}$
 of $X_{0,\overline{\eta}}$ in $X_{\overline{\eta}}$ is proper over $S$.
By Lemma \ref{lem:verify-assump} i), ii), we can apply Theorem \ref{thm:LTF-open-adic} to 
 $X_{0,\overline{\eta}}\hooklongrightarrow \overline{X_{0,\overline{\eta}}}$. 
 Together with Proposition  \ref{prop:trace-identity} and Lemma \ref{lem:verify-assump} iii), we can conclude
 \[
  \Tr\bigl(f_{\overline{\eta}}^*;R\Gamma_c(X_{\overline{\eta}},\Lambda)\bigr)=
 \Tr\bigl(f^*_{\overline{\eta}}\vert_{X_{0,\overline{\eta}}};R\Gamma_c(X_{0,\overline{\eta}},\Lambda)\bigr)
 =\#\Fix (f_{\overline{\eta}}\vert_{X_{0,\overline{\eta}}})=\#\Fix f_{\overline{\eta}}.
 \]
\end{prf}

\begin{rem}
 At least when the characteristic of $k$ is $0$,
 we can deduce from Theorem \ref{thm:LTF-formal} the analogous result for $\ell$-adic coefficient 
 simply by taking projective limit (\cf \cite[proof of Corollary 4.40]{formalnearby}).
\end{rem}

Next we discuss the existence of a system of neighborhoods in
Assumption \ref{assump:morphism}. Let $f\colon \mathcal{X}\yrightarrow{\cong}\mathcal{X}$
and $\{\mathcal{Y}_\alpha\}_{\alpha\in \mathcal{T}}$ be as in the beginning of this section, and
$f\colon \mathcal{T}\yrightarrow{\cong} \mathcal{T}$ a bijection (we do not need to take the order on $\mathcal{T}$ 
into account).
We want to find a system of closed constructible subsets $\{Y_\alpha(n)\}_{n\ge 1}$
satisfying i), ii), iii) in Assumption \ref{assump:morphism}.
To construct it, we introduce a ``tubular neighborhood'' of a closed formal subscheme $\mathcal{Y}$ of $\mathcal{X}$.

\begin{defn}\label{defn:tubular}
 Let $\mathcal{Y}$ be a closed formal subscheme of $\mathcal{X}$
 and $\mathcal{I}$ be an ideal of definition of $\mathcal{X}$.
 We will define the subset $Y(\mathcal{I})$ of $X=t(\mathcal{X})_a$ as follows.
 First assume that $\mathcal{X}=\Spf A$ is affine. Then $\mathcal{Y}$ is defined by an ideal $J$ of $A$.
 Put $I=\Gamma(\mathcal{X},\mathcal{I})$ and
 \[
 Y(\mathcal{I})=\bigl\{x\in X\bigm| \max_{f\in J}\lvert f(x)\rvert<\max_{g\in I}\lvert g(x)\rvert\bigr\}.
 \]
 Note that $\max_{f\in J}\lvert f(x)\rvert=\max_{1\le i\le m}\lvert f_i(x)\rvert$ for every system of generators
 $f_1,\ldots,f_m$ of $J$; in particular $\max_{f\in J}\lvert f(x)\rvert$ exists.
 Similar for $\max_{g\in I}\lvert g(x)\rvert$.

 Obviously we can globalize the construction by patching, and get $Y(\mathcal{I})$
 in the general case.
\end{defn}

\begin{prop}\label{prop:tubular-constr}
 The subset $Y(\mathcal{I})$ is closed constructible in $X$.
\end{prop}

\begin{prf}
 We may assume that $\mathcal{X}=\Spf A$ is affine. Let $J\subset A$ be the defining ideal of $\mathcal{Y}$ and
 put $I=\Gamma(\mathcal{X},\mathcal{I})$. Take a system of generators $f_1,\ldots,f_m$ (resp.\ $g_1,\ldots,g_n$)
 of $J$ (resp.\ $I$). Then, noting that $X=\{x\in t(\mathcal{X})\mid \max_{1\le j\le n}\lvert g_j(x)\rvert\neq 0\}$,
 we have
 \[
  X\setminus Y(\mathcal{I})=X\cap \bigcup_{1\le i\le m}R\Bigl(\frac{g_1,\ldots,g_n}{f_i}\Bigr),
 \]
 where $R(-)$ denotes a rational subset of $t(\mathcal{X})=\Spa(A,A)$.
 Since every rational subset is quasi-compact and open,
 $X\setminus Y(\mathcal{I})$ is a quasi-compact open subset of $X$. 
 This completes the proof.
\end{prf}

The following two lemmas are clear from the definition:

\begin{lem}\label{lem:tubular-fundamental}
 For an ideal of definition $\mathcal{I}$ of $\mathcal{X}$, we have the following:
 \begin{enumerate}
  \item $Y(\mathcal{I}^{n+1})\subset Y(\mathcal{I}^n)$
	for every $n\ge 1$.
  \item $Y:=t(\mathcal{Y})_a=\bigcap_{n\ge 1}Y(\mathcal{I}^n)$.
 \end{enumerate}
\end{lem}

\begin{lem}\label{lem:tubular-functoriality}
 Let $f\colon \mathcal{X}'\longrightarrow\mathcal{X}$ be an adic morphism over $\mathcal{S}$ and 
 put $\mathcal{Y}'=\mathcal{Y}\times_{\mathcal{X}}\mathcal{X}'$.
 For an ideal of definition $\mathcal{I}$ of $\mathcal{X}$, 
 put $\mathcal{I}'=(f^{-1}\mathcal{I})\mathcal{O}_{\mathcal{X'}}$.
 Then $f\colon X':=t(\mathcal{X'})_a\longrightarrow X$ induces a map from $Y'(\mathcal{I}')$ to $Y(\mathcal{I})$.
\end{lem}

Now we can give fairly simple conditions for existence of a system of neighborhoods
in Assumption \ref{assump:morphism}.

\begin{prop}\label{prop:suff-cond-isom}
 Assume that the isomorphism $f\colon \mathcal{X}\yrightarrow{\cong}\mathcal{X}$
 induces an isomorphism of formal schemes $\mathcal{Y}_\alpha\yrightarrow{\cong} \mathcal{Y}_{f(\alpha)}$
 for every $\alpha\in \mathcal{T}$. Then there exists a system of closed constructible subsets
 $\{Y_\alpha(n)\}_{n\ge 1}$ of $X$ satisfying i), ii), iii)
 in Assumption \ref{assump:morphism}.
\end{prop}

\begin{prf}
 Let $\mathcal{I}$ be the maximal ideal of definition of $\mathcal{X}$ (it exists since $\mathcal{X}$ is noetherian)
 and put $Y_\alpha(n)=Y_\alpha(\mathcal{I}^n)$ for each $\alpha\in\mathcal{T}$ and $n\ge 1$.
 By Proposition \ref{prop:tubular-constr}, $Y_\alpha(n)$ is a closed constructible 
 subset of $X$.
 Moreover, by Lemma \ref{lem:tubular-fundamental},
 $\{Y_\alpha(n)\}_{n\ge 1}$ satisfies the conditions i), ii) in Assumption \ref{assump:morphism}.
 Finally, the condition iii) follows from Lemma \ref{lem:tubular-functoriality},
 since $\mathcal{I}$ is preserved by the isomorphism $f$.
\end{prf}

\begin{prop}\label{prop:suff-cond-top-fin-order}
 Assume that the isomorphism $f\colon \mathcal{X}\yrightarrow{\cong}\mathcal{X}$
 induces a set-theoretic bijection $Y_\alpha\yrightarrow{\cong}Y_{f(\alpha)}$ for every $\alpha\in \mathcal{T}$.
 Assume moreover that for every ideal of definition $\mathcal{I}$ of $\mathcal{X}$, there exists an integer $N\ge 1$
 such that $f^N\equiv \id \pmod{\mathcal{I}}$. 
 Then there exists a system of closed constructible subsets
 $\{Y_\alpha(n)\}_{n\ge 1}$ of $X$ satisfying i), ii), iii)
 in Assumption \ref{assump:morphism}.
\end{prop}

First we will show:

\begin{lem}
 Assume that the isomorphism $f\colon \mathcal{X}\yrightarrow{\cong}\mathcal{X}$ satisfies the condition
 in Proposition \ref{prop:suff-cond-top-fin-order}. Then, for every constructible subset $V$ of $X$,
 there exists an integer $N\ge 1$ such that $f^N(V)=V$.
\end{lem}

\begin{prf}
 Replacing $f$ by its power if necessary, we may assume that $f$ induces the identity on the underlying space of 
 $\mathcal{X}$.
 Therefore we may assume that $\mathcal{X}$ is affine.
 Moreover, we can reduce to the case where $V$ is a rational subset of $t(\mathcal{X})$.
 Now the lemma is clear from \cite[Lemma 3.10]{MR1207303}.
\end{prf}

\begin{prf}[of Proposition \ref{prop:suff-cond-top-fin-order}]
 Take an ideal of definition $\mathcal{I}$ of $\mathcal{X}$.

 Let us decompose $\mathcal{T}$ into $f$-orbits $\mathcal{T}_1\amalg\cdots\amalg\mathcal{T}_m$.
 Fix an element $\alpha_i\in \mathcal{T}_i$ for each $1\le i\le m$; then 
 $\mathcal{T}_i=\{f^j(\alpha_i)\mid j\in \Z\}$.
 For $\alpha\in\mathcal{T}_i$ and an integer $n\ge 1$, put
 $Y_\alpha(n)=\bigcap_{j\in\Z, f^j(\alpha_i)=\alpha}f^j(Y_{\alpha_i}(\mathcal{I}^n))$.
 By the previous lemma, this intersection is essentially finite.
 Therefore, $Y_\alpha(n)$ is a closed constructible subset of $X$.
 By the assumption $f(Y_\alpha)=Y_{f(\alpha)}$ and Lemma \ref{lem:tubular-fundamental}, 
 $\{Y_\alpha(n)\}_{n\ge 1}$ satisfies the condition i), ii) in Assumption \ref{assump:morphism}.
 On the other hand, $f(Y_\alpha(n))=Y_{f(\alpha)}(n)$ is
 clear from the construction. Now the proof is complete.
\end{prf}

We finish this section by two examples of Rapoport-Zink spaces.

\begin{exa}\label{exa:Lubin-Tate}
 Let $\mathcal{O}$ be a complete discrete valuation ring with finite residue field $\F_q$.
 Denote the completion of the strict henselization of $\mathcal{O}$ by $\breve{\mathcal{O}}$ and
 the fraction field of $\breve{\mathcal{O}}$ by $\breve{F}$.
 Fix an integer $d\ge 1$ and denote by $\mathbb{X}$ a formal $\mathcal{O}$-module over $\overline{\mathbb{F}}_q$
 with $\mathcal{O}$-height $d$ (such $\mathbb{X}$ is unique up to isomorphism).
 For an integer $m\ge 0$, $\mathcal{X}_m$ denotes the universal deformation space over 
 $\breve{\mathcal{O}}$ of $\mathbb{X}$ with Drinfeld $m$-level structures.
 For the precise definition, see \cite[\S 2.1]{MR2383890} for example.
 Recall that $\mathcal{X}_0$ is isomorphic to $\Spf \breve{\mathcal{O}}[[T_1,\ldots,T_{d-1}]]$ and
 the natural morphism $\mathcal{X}_m\longrightarrow \mathcal{X}_0$ is finite. In particular,
 $\mathcal{X}_m$ is special over $\Spf \breve{\mathcal{O}}$ and its generic fiber
 $X_m=t(\mathcal{X}_m)_{\eta}$ is partially proper over $\Spa(\breve{F},\breve{\mathcal{O}})$.
 Moreover, it is known that the morphism $X_m\longrightarrow X_0$ induced on the generic fibers is \'etale.
 Therefore $X_m$ is smooth over $\Spa(\breve{F},\breve{\mathcal{O}})$.

 More generally, we can associate to a compact open subgroup $K$ of $K_0=\mathrm{GL}_d(\mathcal{O})$
  the formal scheme $\mathcal{X}_K$ (\cf \cite[\S 2.2]{MR2383890}).
 Put $K_m=\Ker(\mathrm{GL}_d(\mathcal{O})\to \mathrm{GL}_d(\mathcal{O}/\mathfrak{m}^m))$ for an integer $m\ge 1$, 
 where $\mathfrak{m}$ denotes the maximal ideal of $\mathcal{O}$.
 Take $m\ge 0$ such that $K_m\subset K$; then $\mathcal{X}_K$ is defined as the quotient in the sense of invariant theory
 of the action of the finite group $K/K_m$ on $\mathcal{X}_m$ (the action of $K\subset K_0$ on $\mathcal{X}_m$ is
 given via the Drinfeld level structures). 
 It is easy to see that $\mathcal{X}_K$ is special over $\Spf \breve{\mathcal{O}}$ and
 its generic fiber $X_K=t(\mathcal{X}_K)_\eta$ is partially proper and smooth over $\Spa(\breve{F},\breve{\mathcal{O}})$.

 Let $D$ be the central division algebra over $F$ with invariant $1/d$.
 The formal scheme $\mathcal{X}_K$ is endowed with a right action of the subgroup of $\mathrm{GL}_d(F)\times D^\times$
 consisting of elements $(g,h)$ such that $v_F(\det g)+v_F(\Nrd h)=0$ and $gKg^{-1}=K$,
 where $v_F$ denotes the normalized valuation of $F$. 
 We would like to explain that we can use Theorem \ref{thm:LTF-formal} to calculate the trace
 $\Tr((g,h)^*;R\Gamma_c(X_{K,\overline{\eta}},\Lambda))$, under the assumption that $gK$ consists of
 regular elliptic elements of $\mathrm{GL}_d(F)$.

 For an integer $m\ge 1$, let $\mathcal{S}_m$ be the ordered set of $\mathcal{O}/\mathfrak{m}^m$-submodules of 
 $(\mathfrak{m}^{-m}/\mathcal{O})^d$ which are direct summands.
 Put $\mathcal{S}_m'=\mathcal{S}_m\setminus \{0,(\mathfrak{m}^{-m}/\mathcal{O})^d\}$.
 For each $I\in \mathcal{S}_m$, we can construct the closed formal subscheme $\mathcal{Y}_I$ of $\mathcal{X}_{m,s}$
 (\cf \cite[Definition 4.1]{non-cusp});
 roughly speaking, it is the locus where the universal Drinfeld $m$-level structure vanishes on $I$.
 More generally, for an compact open subgroup $K$ of $K_0$, we put $\mathcal{S}_K=(K/K_m)\backslash\mathcal{S}_m$ and 
 $\mathcal{S}'_K=(K/K_m)\backslash\mathcal{S}'_m$, where $m\ge 1$ is an integer with $K_m\subset K$.
 We endow them with the induced partial orders. For $I\in \mathcal{S}_K$,
 we can also define the closed formal subscheme $\mathcal{Y}_I$ of $\mathcal{X}_{K,s}$
 so that $Y_{(I)}=Y_I\setminus \bigcup_{I'\in\mathcal{S}_K, I'>I}Y_{I'}$, where $Y_I$ denotes $t(\mathcal{Y}_I)_a$,
 coincides with the boundary subset $\partial_I M_K$ in \cite[Paragraph 3.1.1]{MR2383890};
 in fact, we can take $\mathcal{Y}_I$ as the closed formal subscheme defined by $\mathfrak{p}_K$
 in \cite[Proposition 3.1.3 (i)]{MR2383890}
 (note that the partial order on $\mathcal{S}_K$ in \cite{MR2383890} is the inverse of ours).
 However, the author could not find any natural moduli interpretation of $\mathcal{Y}_I$.

 Now we can easily see that the family of closed formal subschemes $\{\mathcal{Y}_I\}_{I\in\mathcal{S}'_m}$ satisfies
 Assumption \ref{assump:formal-sch}.
 On the other hand, we can define the action of $g$ on $\mathcal{S}'_K$ as in \cite[p.~914]{MR2383890}.
 Since we are assuming that $gK$ consists of regular elliptic elements,
 we have $I\neq g^{-1}I$ for every $I\in\mathcal{S}'_K$.
 Moreover, it is easy to see that the action of $(g,h)$ maps $Y_I$ onto $Y_{g^{-1}I}$ at least set-theoretically
 (\cf \cite[Lemma 3.2.2 (ii)]{MR2383890}). Since $(g,h)$ is elliptic,
 the action of $(g,h)$ on $\mathcal{X}_K$ satisfies the condition in Proposition \ref{prop:suff-cond-top-fin-order}
 by the first part of the proof of \cite[Proposition 3.2.4]{MR2383890}.
 Therefore, by Proposition \ref{prop:suff-cond-top-fin-order},
 the action of $(g,h)$ on $\mathcal{X}_K$ satisfies Assumption \ref{assump:morphism}.
 Finally, it is known that $\mathcal{X}_K$ is algebraizable (\cite[Theorem 2.3.1]{MR2383890}).
 Hence all the assumptions in Theorem \ref{thm:LTF-formal} are satisfied and we obtain
 \[
  \Tr\bigl((g,h)^*;R\Gamma_c(X_{K,\overline{\eta}},\Lambda)\bigr)=\#\Fix(g,h).
 \]
 This recovers a result \cite[Theorem 3.3.1]{MR2383890} of Strauch.
 Recall that the right hand side has been calculated in \cite[Theorem 2.6.8]{MR2383890},
 and as a consequence of the trace formula above, we can get a purely local proof of the fact that the
 $\ell$-adic cohomology of the Lubin-Tate
 tower $(X_m)_{m\ge 0}$ realizes the local Jacquet-Langlands correspondence 
 (\cite[Theorem 4.1.3]{MR2383890}; see also \cite{LT-LTF}).

 The advantage of our proof is that it does not require algebraization of the action.
 The proof of \cite[Theorem 3.3.1]{MR2383890} uses careful approximation of the action of $(g,h)$
 by an algebraizable morphism (\cf \cite[Proposition 3.2.4 (ii), \S 5.2]{MR2383890}),
 and it seems difficult to extend that method to the non-affine case.
\end{exa}

\begin{exa}
 Let $p$ be a prime. For a compact open subgroup $K^p$ of $\mathrm{GSp}_4(\A^{\infty,p})$ and an integer $m\ge 0$,
 let $\Sh_{m,K^p}$ be the Shimura variety over $\Z_{p^\infty}=W(\overline{\F}_p)$
 introduced in \cite[\S 4]{RZ-GSp4}; namely, it is the moduli space parameterizing polarized abelian surfaces
 with $K^p$-level structures outside $p$ and Drinfeld $m$-level structures at $p$.
 Let $\mathcal{S}_m$ be the ordered set of direct summands of $(\Z/p^m\Z)^4$ whose ranks are greater than $1$
 and put $\mathcal{S}'_m=\mathcal{S}_m\setminus \{(\Z/p^m\Z)^4\}$.
 Fix a perfect alternating bilinear form on $\Z_p^4$ and
 denote the subset of $\mathcal{S}_m$ (resp.\ $\mathcal{S}'_m$) consisting of coisotropic direct summands
 by $\mathcal{S}_m^{\mathrm{coi}}$ (resp.\ $\mathcal{S}_m'^{\mathrm{coi}}$). 
 Then, as in Example \ref{exa:Lubin-Tate}, we can define the closed subscheme 
 $\overline{\Sh}_{m,K^p,[I]}$ of $\overline{\Sh}_{m,K^p}=\Sh_{m,K^p}\otimes_{\Z_{p^\infty}}\overline{\F}_p$ for
 $I\in\mathcal{S}_m$ (\cf \cite[Definition 5.1]{RZ-GSp4}).
 It is known that every $\overline{\F}_p$-rational point of $Y_{m,K^p}:=\overline{\Sh}_{m,K^p,[(\Z/p^m\Z)^4]}$
 corresponds to a supersingular abelian variety (note that the definition of $Y_{m,K^p}$ here is different from 
 that in \cite{RZ-GSp4}, but they coincide up to nilpotent elements; see \cite[Lemma 5.3 iii)]{RZ-GSp4}).

 Let $\Sh'_{m,K^p}$ be the closed subscheme of $\Sh_{m,K^p}$ defined by the quasi-coherent ideal 
 of $\mathcal{O}_{\Sh_{m,K^p}}$ consisting of elements killed by $p^l$ for some integer $l\ge 0$.
 Put $\overline{\Sh}^{\,\prime}_{m,K^p}=\overline{\Sh}_{m,K^p}\times_{\Sh_{m,K^p}}\Sh'_{m,K^p}$,
 $\overline{\Sh}^{\,\prime}_{m,K^p,[I]}=\overline{\Sh}_{m,K^p,[I]}\times_{\Sh_{m,K^p}}\Sh'_{m,K^p}$ and
 $Y'_{m,K^p}=Y_{m,K^p}\times_{\Sh_{m,K^p}}\Sh'_{m,K^p}$. Then, by \cite[Lemma 5.3, proof of Proposition 5.7]{RZ-GSp4}
 the family of closed subschemes $\{\overline{\Sh}^{\,\prime}_{m,K^p,[I]}\}_{I\in\mathcal{S}^{\mathrm{coi}}_m}$
 satisfies the conditions in Example \ref{exa:assump-scheme}.
 Let us denote the completion of $\Sh_{m,K^p}$ (resp.\ $\Sh'_{m,K^p}$) along $Y_{m,K^p}$ (resp.\ $Y'_{m,K^p}$)
 by $(\Sh_{m,K^p})^{\wedge}_{/Y_{m,K^p}}$ (resp.\ $(\Sh'_{m,K^p})^{\wedge}_{/Y'_{m,K^p}}$).

 Let $\M$ be the Rapoport-Zink space for $\mathrm{GSp}(4)$ considered in \cite{RZ-GSp4} and
 $\M_m$ be as in \cite[\S 3.2]{RZ-GSp4}. The formal scheme $\M_m$ is endowed with the action of 
 $\mathrm{GSp}_4(\Z_p)\times J(\Q_p)$, where $J$ is an inner form of $\mathrm{GSp}_4$.
 By the $p$-adic uniformization theorem of Rapoport-Zink, 
 there is a $\mathrm{GSp}_4(\Z_p)$-equivariant morphism
 $\theta_{m,K^p}\colon \M_m\longrightarrow (\Sh_{m,K^p})^{\wedge}_{/Y_{m,K^p}}$,
 which induces an open and closed immersion 
 $\M_m/\Gamma\hooklongrightarrow (\Sh_{m,K^p})^{\wedge}_{/Y_{m,K^p}}$ for some
 discrete subgroup $\Gamma$ of $J(\Q_p)$ (\cite[Theorem 4.2]{RZ-GSp4}).
 Put $\M_{m,[I]}=\M_m\times_{\Sh_{m,K^p}}\overline{\Sh}_{m,K^p,[I]}$,
 $\M^{\,\prime}_m=\M_m\times_{\Sh_{m,K^p}}\Sh'_{m,K^p}$ and 
 $\M^{\,\prime}_{m,[I]}=\M_m\times_{\Sh_{m,K^p}}\overline{\Sh}^{\,\prime}_{m,K^p,[I]}$.
 As mentioned in \cite[proof of Proposition 5.10]{RZ-GSp4}, $\M_{m,[I]}$ is preserved by the action of $J$ on $\M_m$.
 On the other hand, it is easy to see that the defining ideal of $\M^{\,\prime}_m$ in $\M_m$ consists of
 elements of $\mathcal{O}_{\!\!\M_m}$ which is killed by $p^l$ for some integer $l>0$. Therefore the actions 
 of $J$ on $\M^{\,\prime}_m$ and $\M^{\,\prime}_{m,[I]}$ are naturally induced.
 Moreover, we have an open and closed immersion $\M^{\,\prime}_m/\Gamma\hooklongrightarrow (\Sh'_{m,K^p})^{\wedge}_{/Y'_{m,K^p}}$. It is also easy to observe that $t(\M_m^{\,\prime})_\eta$ coincides with $t(\M_m)_\eta$.

 Now it is easy to see that the formal scheme $\M^{\,\prime}_m/\Gamma$, a family of closed formal subscheme
 $\{\M^{\,\prime}_{m,[I]}/\Gamma\}_{I\in\mathcal{S}_m'^{\mathrm{coi}}}$ and the action of $(g,j)\in \mathrm{GSp}_4(\Z_p)\times J(\Q_p)$ where $gK_m$ consists of regular elliptic elements of $\mathrm{GSp}_4(\Q_p)$ and $j$ normalizes $\Gamma$
 satisfy all the assumptions in Theorem \ref{thm:LTF-formal}.
 Indeed, Assumption \ref{assump:formal-sch} has already been observed.
 Assumption \ref{assump:morphism} follows from \cite[Proposition 5.15]{RZ-GSp4}, Proposition \ref{prop:suff-cond-isom} 
 and the fact that if $gK_m$ consists of regular elliptic then $g^{-1}I\neq I$
 for every $I\in\mathcal{S}_m'^{\mathrm{coi}}$
 (otherwise $gK_m$ intersects a proper parabolic subgroup).
 Since every irreducible component of $(\M_m)_{\mathrm{red}}$ is
 projective over $\overline{\F}_p$ (\cite[Proposition 2.32]{MR1393439}), $Y_{m,K^p}$ and $Y'_{m,K^p}$ are proper
 over $\overline{\F}_p$ and thus $\M_m^{\,\prime}/\Gamma$ is partially proper over $\Spa(\Q_{p^\infty},\Z_{p^\infty})$.
 Hence we get the formula
 \[
  \Tr\bigl((g,j)^*;R\Gamma_c(t(\M_m)_{\overline{\eta}}/\Gamma,\Lambda)\bigr)=\#\Fix(g,j).
 \]
 The right hand side can be calculated in the similar way as in \cite[\S 2.6]{MR2383890}.
 The author plans to work this out in a forthcoming paper.
\end{exa}

\section{Lefschetz trace formula for contracting morphisms}\label{sec:LTF-contr}
\subsection{Statement}\label{subsec:contracting-statement}
In this section, we generalize Fujiwara's trace formula for contracting morphisms 
(\cite[Theorem 3.2.4]{MR1431137}) to rigid spaces
which are not necessarily algebraizable. Let $X$ be a purely $d$-dimensional adic space which is
proper and smooth over $S$. For a closed adic subspace $Y$ of $X$ and $\varepsilon\in \lvert k^\times\rvert\subset \R$,
we can construct the open tubular neighborhood $Y^{\circ}(\varepsilon)$ and the closed tubular neighborhood
$Y(\varepsilon)$ (\cf \cite[\S 2.6]{MR1620114}, in which $Y^{\circ}(\varepsilon)$ is denoted by
$T(\varepsilon)$ and $Y(\varepsilon)$ by $S(\varepsilon)$).
If $Y$ is defined by $f_1,\ldots,f_m\in \Gamma(X,\mathcal{O}_X)$, then 
$Y^{\circ}(\varepsilon)$ (resp. $Y(\varepsilon)$) is given by
$\{x\in X\mid \lvert f_i(x)\rvert\le \varepsilon\}$ (resp.\ $\{x\in X\mid \lvert f_i(x)\rvert<\varepsilon\}$).
Note that $Y^{\circ}(\varepsilon)$ and $Y(\varepsilon)$ are constructible and
$Y=\bigcap_{\varepsilon\in\lvert k^\times\rvert}Y^\circ(\varepsilon)=\bigcap_{\varepsilon\in\lvert k^\times\rvert}Y(\varepsilon)$.

Let $f\colon X\longrightarrow X$ be an $S$-morphism. We will use the notation in Example \ref{exa:self-mor}.
We denote the set of connected components of $\Fix f$ by $\pi_0(\Fix f)$.
It is a finite set since $H^0(\Fix f,\Lambda)$ is a finitely generated $\Lambda$-module. Therefore
every element of $\pi_0(\Fix f)$ is open and closed in $\Fix f$.

\begin{defn}
 Let $D$ be a connected component of $\Fix f$.
 We say that $f$ is contracting near $D$ if there exists a strictly decreasing sequence
 $(\varepsilon_n)_{n\ge 0}$ in $\lvert k^\times\rvert$ converging to $0$ such that 
 $f(D^\circ(\varepsilon_n))\subset D^\circ(\varepsilon_{n+1})$ for every $n\ge 0$.

 If $f$ is contracting near every connected component of $\Fix f$, we say that $f$ is contracting near
 its fixed points.
\end{defn}

The goal of this section is the following theorem:

\begin{thm}\label{thm:LTF-contr}
 Assume that the characteristic of $k$ is equal to $0$ and $f$ is contracting near $D\in \pi_0(\Fix f)$.
 Then we have
 \[
 \#\Fix_D f=\chi(D,\Lambda),
 \]
 where we put $\chi(D,\Lambda)=\Tr(\id; R\Gamma(D,\Lambda))$.
\end{thm}

This theorem will be proved in \S \ref{subsec:proof-contracting}. The following corollary is immediate from
Theorem \ref{thm:LTF-contr} and Theorem \ref{thm:LTF-open-adic}:

\begin{cor}
 Assume that the characteristic of $k$ is equal to $0$ and $f$ is contracting near its fixed points.
 Then we have
 \[
 \Tr\bigl(f^*;R\Gamma(X,\Lambda)\bigr)=\sum_{D\in \pi_0(\Fix f)}\chi(D,\Lambda).
 \]
 In particular, if every fixed point of $f$ is isolated, the right hand side is equal to the number of the fixed points.
\end{cor}

\subsection{Lefschetz trace formula for proper pseudo-adic spaces}
In order to prove Theorem \ref{thm:LTF-contr}, we need a variant of Proposition \ref{prop:LTF-unlocalized}
for a pseudo-adic space which is proper but not necessarily smooth over $S$.
In this subsection, let $X$ be a finite-dimensional pseudo-adic space which is proper over $S$.
We assume that $H^i(X,\Lambda)$ is a finitely generated $\Lambda$-module for every integer $i$;
then $R\Gamma(X,\Lambda)$ is a perfect $\Lambda$-complex (Corollary \ref{cor:perfect-cpx}).

Let $f\colon X\longrightarrow X$ be an $S$-morphism.
Since $R\Gamma(X,\Lambda)$ is a perfect complex, the trace
$\Tr(f^*;R\Gamma(X,\Lambda))$ makes sense. Our purpose is to express this trace by a cohomology class
analogous to $[\gamma]$ in Proposition \ref{prop:LTF-unlocalized}.

First we will construct an analogue of $\cl(\gamma)$. 

\begin{defn}
 Since $\pr_2\circ \gamma_f=\id$, we have $H^0_{\Gamma_f}(X\times_SX,\pr_2^!\Lambda)\cong H^0(X,\Lambda)$.
 We denote by $\cl(f)$ the element of $H^0_{\Gamma_f}(X\times_SX,\pr_2^!\Lambda)$ that corresponds to
 $1\in H^0(X,\Lambda)$ under the isomorphism above.
\end{defn}

\begin{rem}\label{rem:cl-cl}
 If $X$ is purely $d$-dimensional and smooth over $S$,
 then $H^0_{\Gamma_f}(X\times_SX,\pr_2^!\Lambda)\cong H^{2d}_{\Gamma_f}(X\times_SX,\Lambda(d))$ and
 by this isomorphism $\cl(f)$ corresponds to $\cl(\gamma_f)$ in \S 2.
\end{rem}

The following is an analogue of Proposition \ref{prop:corr-cl}, whose proof is similar:

\begin{prop}\label{prop:pull-cl}
 The map $f^*\colon R\Gamma(X,\Lambda)\longrightarrow R\Gamma(X,\Lambda)$ coincides with the composite of
 \begin{align*}
  R\Gamma(X,\Lambda)&\yrightarrow{\pr_1^*}R\Gamma(X\times_SX,\Lambda)\yrightarrow{\cup\cl(f)}
 R\Gamma(X\times_SX,\pr_2^!\Lambda)=R\Gamma(X,R\pr_{2!}\pr_2^!\Lambda)\\
 &\yrightarrow{\adj}R\Gamma(X,\Lambda).
 \end{align*}
\end{prop}

Next we will establish an analogue of Corollary \ref{cor:Kunneth-square-shape}. 
Denote the structure morphism of $X$ by $a\colon X\longrightarrow S$ and put $K_X=a^!\Lambda$.
Let $\tau\colon K_X\Lboxtimes \Lambda\longrightarrow \pr_2^!\Lambda$ be the base change map
$K_X\Lboxtimes \Lambda=\pr_1^*a^!\Lambda\longrightarrow \pr_2^!a^*\Lambda=\pr_2^!\Lambda$.

\begin{prop}\label{prop:Kunneth-K}
 The map $\tau$ above induces an isomorphism
 \[
  R\Gamma(X\times_SX,K_X\Lboxtimes \Lambda)\yrightarrow{\cong} R\Gamma(X\times_SX,\pr_2^!\Lambda).
 \]
\end{prop}

\begin{prf}
 By Proposition \ref{prop:Kunneth-const}, the K\"unneth morphism
 $R\Gamma(X,K_X)\Lotimes R\Gamma(X,\Lambda)\longrightarrow R\Gamma(X\times_SX,K_X\Lboxtimes \Lambda)$ is
 an isomorphism.
 On the other hand, we have isomorphisms
 \begin{align*}
  &R\Gamma(X\times_SX,\pr_2^!\Lambda)=R\Gamma(X,R\pr_{1*}\pr_2^!\Lambda)
 \cong R\Gamma(X,a^!Ra_*\Lambda)=R\Gamma\bigl(X,a^!R\Gamma(X,\Lambda)\bigr)\\
  &\qquad \underset{(1)}{\cong} R\Gamma\bigl(X,a^!\Lambda\Lotimes R\Gamma(X,\Lambda)_X\bigr)
  \underset{(2)}{\cong} R\Gamma(X,a^!\Lambda)\Lotimes R\Gamma(X,\Lambda)\\
  &\qquad =R\Gamma(X,K_X)\Lotimes R\Gamma(X,\Lambda).
 \end{align*}
 For (1), note that the natural map $a^!\Lambda\Lotimes R\Gamma(X,\Lambda)_X\longrightarrow a^!R\Gamma(X,\Lambda)$
 is an isomorphism. This is an easy consequence of the fact that $R\Gamma(X,\Lambda)$
 is a perfect $\Lambda$-complex. The isomorphy of (2) is due to Lemma \ref{lem:constant}.

 It is easy to see that this isomorphism fits into the following diagram:
 \[
  \xymatrix{%
 R\Gamma(X,K_X)\Lotimes R\Gamma(X,\Lambda)\ar@<-3pt>@{=}[r]\ar[d]^-{\cong}&R\Gamma(X,K_X)\Lotimes R\Gamma(X,\Lambda)\ar[d]^-{\cong}\\
 R\Gamma(X\times_SX,K_X\Lboxtimes \Lambda)\ar[r]^-{\tau}& R\Gamma(X\times_SX,\pr_2^!\Lambda)\lefteqn{.}
 }
 \]
 This completes the proof.
\end{prf}

\begin{defn}
 Let $[f]$ be the element of $H^0(X\times_SX,K_X\Lboxtimes\Lambda)$ that corresponds to $\cl(f)$
 by the isomorphism in Proposition \ref{prop:Kunneth-K}.
\end{defn}

\begin{prop}\label{prop:LTF-nonsmooth}
 We have $\Tr(f^*;R\Gamma(X,\Lambda))=\Adj_X(\delta^*[f])$, where $\Adj_X$ denotes the natural adjunction map
 $H^0(X,K_X)\longrightarrow \Lambda$.
\end{prop}

\begin{prf}
 By Proposition \ref{prop:pull-cl}, it suffices to show the commutativity of the lower part of the following diagram:
 \[
 \xymatrix{%
 R\Gamma(X,K_X)\Lotimes R\Gamma(X,\Lambda)\ar[d]^-{\cong}
 \ar@<-2pt>[r]^-{\cong}& R\Hom\bigl(R\Gamma(X,\Lambda),\Lambda\bigr)\Lotimes R\Gamma(X,\Lambda)\ar[dd]\\
 R\Gamma(X\times_SX,K_X\Lboxtimes \Lambda)\ar[d]\ar@/_6pc/[dd]_-{\delta^*}\\
 R\Gamma(X\times_SX,\pr_2^!\Lambda)\ar[r]^-{(*)}& R\Hom\bigl(R\Gamma(X,\Lambda),
 R\Gamma(X,\Lambda)\bigr)\ar[d]^-{\Tr}\\
 R\Gamma(X,K_X)\ar[r]^-{\Adj_X}& \Lambda\lefteqn{,}
 }
\]
 where $(*)$ is given by
 \begin{align*}
  &R\Gamma(X,\Lambda)\Lotimes R\Gamma(X\times_SX,\pr_2^!\Lambda)\yrightarrow{\pr_1^*\cup \id}R\Gamma(X\times_SX,\pr_2^!\Lambda)\\
  &\qquad=R\Gamma(X,R\pr_{2!}\pr_2^!\Lambda)\yrightarrow{\adj}R\Gamma(X,\Lambda).
 \end{align*}
 As in the proof of Proposition \ref{prop:LTF-unlocalized}, we can prove the commutativities of the upper part
 and the outer part of the diagram above.
\end{prf}

We can apply the technique in \S \ref{subsec:localization} to calculate $\Adj_X(\delta^*[f])$:

\begin{lem}\label{lem:calc-Adj}
 Let $U$ be an open adic subspace of $X$ which is purely $d$-dimensional and smooth over $S$, 
 and $Y$ an closed constructible subset of $X$ contained in $U$.
 Assume that $f(X)\subset Y$. Then we have $\Adj_X(\delta^*[f])=\#\Fix f\vert_U$.
\end{lem}

\begin{prf}
 First note that $\Fix f\vert_U=\Fix f$ since $\Fix f\subset Y$. 
 Therefore the adic space $\Fix f\vert_U$ is proper over $S$ and $\#\Fix f\vert_U$ makes sense.

 In the same way as Corollary \ref{cor:square-isom-closed-union}, we can deduce from Proposition \ref{prop:Kunneth-K} that
 the map
 $R\Gamma_{Y\times_S X}(X\times_SX,K_X\Lboxtimes\Lambda)\longrightarrow R\Gamma_{Y\times_SX}(X\times_SX,\pr_2^!\Lambda)$ is
 an isomorphism. 

 Consider the following commutative diagram:
 \[
  \xymatrix{%
 H^0_{\Gamma_f}(X\times_SX,\pr_2^!\Lambda)\ar[r]\ar[d]&H^{2d}_{\Gamma_{f\vert_U}}\bigl(U\times_SU,\Lambda(d)\bigr)\ar[d]\\
 H^0_{Y\times_S X}(X\times_SX,\pr_2^!\Lambda)\ar[r]&H^{2d}_{Y\times_S U}\bigl(U\times_SU,\Lambda(d)\bigr)\ar@{=}[d]\\
 H^0_{Y\times_S X}(X\times_SX,K_X\Lboxtimes\Lambda)\ar[r]\ar[d]^-{\delta^*}\ar[u]_{\cong}&H^{2d}_{Y\times_S U}\bigl(U\times_SU,\Lambda(d)\bigr)\ar[d]^-{\delta_U^*}\\
 H^0_Y(X,K_X)\ar[r]^-{\cong}\ar[d]&H^{2d}_Y\bigl(U,\Lambda(d)\bigr)\ar[d]\\
 H^0(X,K_X)\ar[d]^-{\Adj_X}&H^{2d}_c\bigl(U,\Lambda(d)\bigr)\ar[l]\ar[d]^-{\Tr_U}\\
 \Lambda\ar@{=}[r]& \Lambda\lefteqn{.}
 }
 \]
 It is easy to see that the image of $\cl(f)\in H^0_{\Gamma_f}(X\times_SX,\pr_2^!\Lambda)$ under the composite of the arrows in the left column
 is equal to $\Adj_X(\delta^*[f])$. 
 On the other hand, by Remark \ref{rem:cl-cl}, the image of $\cl(f)$ under the top horizontal arrow is $\cl(f\vert_U)$.
 Since the image of $\cl(f\vert_U)$ under the composite of the arrows in the right column is $\#\Fix f\vert_U$,
 we get the lemma.
\end{prf}

\subsection{Proof of Theorem \ref{thm:LTF-contr}}\label{subsec:proof-contracting}
We go back to the notation introduced in \S \ref{subsec:contracting-statement}. 
Let $D$ be a connected component of $\Fix f$ and assume that $f$ is contracting near $D$.
Take a strictly decreasing sequence $(\varepsilon_n)_{n\ge 0}$ in $\lvert k^\times\rvert$ converging to $0$ such that 
$f(D^\circ(\varepsilon_n))\subset D^\circ(\varepsilon_{n+1})$ for every $n\ge 0$.
Take a sequence $(\varepsilon'_n)_{n\ge 0}$ in $\lvert k^\times\rvert$ such that $\varepsilon_n>\varepsilon'_n>\varepsilon_{n+1}$ for every $n\ge 0$.
Then we have $f(D(\varepsilon_n))\subset f(D^\circ(\varepsilon_n))\subset D^\circ(\varepsilon_{n+1})\subset D(\varepsilon'_n)\subset D^\circ(\varepsilon'_n)\subset D(\varepsilon_n)$. 
Fix an integer $n\ge 0$.
The pseudo-adic space $D(\varepsilon_n):=(X,D(\varepsilon_n))$ is proper over $S$ and
$H^i(D(\varepsilon_n),\Lambda)$ is a finitely generated
$\Lambda$-module (\cite[Corollary 2.3]{MR1620114}, \cite[Corollary 5.4]{MR2336836}).
Hence we may apply Proposition \ref{prop:LTF-nonsmooth} to $D(\varepsilon_n)$. 
Moreover, we can also apply Lemma \ref{lem:calc-Adj} to
$D(\varepsilon'_n)\subset D^\circ(\varepsilon'_n)\subset D(\varepsilon_n)$. Summing up, we get the formula
$\Tr(f^*;R\Gamma(D(\varepsilon_n),\Lambda))=\#\Fix (f\vert_{D^\circ(\varepsilon'_n)})$.

On the other hand, for every $D'\in\pi_0(\Fix f)$ distinct from $D$, $D^\circ(\varepsilon'_n)$ does not intersect $D'$;
otherwise $D'$ also intersects $\bigcap_{m\ge 0}f^m(D^\circ(\varepsilon'_n))\subset \bigcap_{m\ge n}D^\circ(\varepsilon_m)=D$.
Thus we have $\Fix (f\vert_{D^\circ(\varepsilon'_n)})=D$ and $\#\Fix (f\vert_{D^\circ(\varepsilon'_n)})=\#\Fix_D (f\vert_{D^\circ(\varepsilon'_n)})=\#\Fix_Df$
(the final equality is due to Proposition \ref{prop:Fix-etale}).

Now assume that the characteristic of $k$ is $0$. Then, by \cite[Theorem 3.6]{MR1620118}, 
the restriction map $H^i(D(\varepsilon_n),\Lambda)\longrightarrow H^i(D,\Lambda)$ is an isomorphism for sufficiently large $n$.
Therefore, by the commutative diagram
 \[
 \xymatrix{%
 R\Gamma\bigl(D(\varepsilon_n),\Lambda\bigr)\ar[r]^-{f^*}\ar[d]^-{\cong}& R\Gamma\bigl(D(\varepsilon_n),\Lambda\bigr)\ar[d]^-{\cong}\\
 R\Gamma(D,\Lambda)\ar[r]^-{\id}& R\Gamma(D,\Lambda)\lefteqn{,}
 }
 \]
 we have $\Tr(f^*;R\Gamma(D(\varepsilon_n),\Lambda))=\chi(D,\Lambda)$ for such $n$.
 Hence we have $\#\Fix_Df=\chi(D,\Lambda)$, as desired.

\def\cftil#1{\ifmmode\setbox7\hbox{$\accent"5E#1$}\else
  \setbox7\hbox{\accent"5E#1}\penalty 10000\relax\fi\raise 1\ht7
  \hbox{\lower1.15ex\hbox to 1\wd7{\hss\accent"7E\hss}}\penalty 10000
  \hskip-1\wd7\penalty 10000\box7} \def\cprime{$'$} \def\cprime{$'$}
\providecommand{\bysame}{\leavevmode\hbox to3em{\hrulefill}\thinspace}
\providecommand{\MR}{\relax\ifhmode\unskip\space\fi MR }
\providecommand{\MRhref}[2]{%
  \href{http://www.ams.org/mathscinet-getitem?mr=#1}{#2}
}
\providecommand{\href}[2]{#2}

\end{document}